\documentclass[12pt]{article}
\usepackage{mathrsfs}
\usepackage{amssymb,a4}
\usepackage{amsmath,amsfonts,amssymb,amsthm}
\newcommand{\al}{\alpha}

\def\wt{{\rm wt}}

\newcommand{\la}{\lambda}

\def\C{{\mathbb C}}

\def\Z{{\mathbb Z}}

\def\1{{\bf 1}}

\def \wt{{\rm wt}}

\def \End{{\rm End}}

\def \<{\langle}
\def \>{\rangle}
\def \a{\alpha}

\def \h{\mathfrak{h}}

\def \w{\omega}
\def \w{\omega}
\def \pf{\noindent {\bf Proof: \,}}
\def\theequation{5.\arabic{equation}}

\renewcommand{\theequation}{\thesection.\arabic{equation}}

\newtheorem{theorem}{Theorem}[section]
\newtheorem{prop}[theorem]{Proposition}
\newtheorem{lem}[theorem]{Lemma}

\newtheorem{remark}[theorem]{Remark}
\theoremstyle{definition}

\begin{document}
\begin{center}
{\Large {\bf  A characterization of vertex operator algebras
$V_{\Z\al}^{+}$: I }}
 \\

\vspace{0.5cm} Chongying Dong\footnote{Supported by NSF grants and a
Faculty research grant from  the University of California at Santa
Cruz.}
\\
Department of Mathematics,  University of California, Santa Cruz, CA 95064 \\
\vspace{.2cm}
\& School of Mathematics,  Sichuan University, Chengdu, 610065 China\\
\vspace{.1 cm} Cuipo Jiang\footnote{Supported  by China NSF grants
10931006,10871125, and the Innovation Program of Shanghai Municipal
Education Commission (11ZZ18).}\\
 Department of Mathematics, Shanghai Jiaotong University, Shanghai 200240 China
\end{center}
\hspace{1cm}

\begin{abstract} A characterization of vertex operator algebra $V_{\Z\alpha}^+$ with $(\alpha,\alpha)/2$
 not being a perfect square is given in terms of dimensions of homogeneous subspaces of small weights.
 This result contributes to the classification of rational vertex operator algebras of central charge $1.$

2000MSC:17B69

\end{abstract}

\section{Introduction}
Classification of rational vertex operator algebras is definitely
one of the most important problems in the theory of vertex operator
algebras and is crucial in classification of rational conformal
field theory. Although there is a lot of progress in the field, one
still has very limited knowledge on the structure theory for vertex
operator algebras. There are two different directions in
classification currently. One direction is the classification of
holomorphic vertex operator algebras which have the simplest
representation theory (see \cite{S}, \cite{DM2}, \cite{LS}). Another
direction is the classification of rational vertex operator algebras
with small central charges. It is established in \cite{DZ} that if
the central charges are less than one, the vertex operator algebra
is an extension of the vertex operator algebra associated to the
discrete series for the Virasoro algebra. But such extensions have
not been constructed and classified except for some special cases.
In the operator algebra setting, classification of local conformal
nets with $c < 1$ has been completed in \cite{KL}. It is natural to
consider classification of rational vertex operator algebras of
central charge 1 next. This has been achieved at the character level
in the physical literature under the assumption that the character
of each irreducible module is a modular function over a congruence
subgroup of the modular group \cite{K}. The classification of
conformal nets of central charge 1 has been given in \cite{X} with
some extra assumption.

It is conjectured that there are three classes rational vertex
operator algebras of central charge 1: (a) vertex operator algebras
$V_L$ associated with positive definite lattices $L$ of rank one,
(b) orbifold vertex operator algebras $V_L^+$ under the automorphism
of $V_L$ induced from the $-1$ isometry of $L$, (c) $V_{\Z\beta}^G$
where $(\beta,\beta)=2$ and $G$ is the finite subgroup of $SO(3)$
isomorphic
 to $A_4, S_4, A_5.$ The lattice vertex operator algebra $V_L$ for any positive definite even
 lattice $L$ has been characterized in \cite{DM1}. This paper gives a characterization of
  $V_L^+$ for $L=\Z\alpha$ with $(\alpha,\alpha)/2$ not being a perfect square. A complete
  characterization of $V_L^+$ for any rank one positive definite even lattice $L$ will be given in a subsequent paper.
  The characterization of $V_{\Z\alpha}^+$ with $(\alpha,\alpha)=4$ has been obtained previously in \cite{ZD} and \cite{DJ}. We should mention that
  the conjecture also requires the effective central charge $\tilde{c}=1$
  (see \cite{ZD}, \cite{DJ}).

There are three major steps in the characterization of $V_L^+.$ Let
$V$ be a rational vertex operator algebra of central charge
 1 satisfying certain conditions. The first step is to show that the Virasoro vector $\omega$ and a weight 4 primary vector $J$
 for the Virasoro algebra generate a subalgebra isomorphic to $M(1)^+$ which is the fixed points of the Heisenberg vertex operator algebra
 of rank one under the $-1$ automorphism. The main idea is to use the fusion rules for the Virasoro
  vertex operator algebra $L(1,0)$ (see \cite{M}, \cite{DJ}) and a result on the $W$-algebra.
The second step is to show that $V$ is a completely reducible
$M(1)^+$-module by using the fusion rules for the vertex operator
algebras $L(1,0)$ and $M(1)^+$ \cite{A1}. The last step is to
establish that the subalgebra $M(1)^+$ and a primary vector
$F\not\in M(1)^+$ of
 minimal weight generate a subalgebra isomorphic to $V_{\Z\alpha}^+$ with $(\alpha,\alpha)=2k$ where $k$ is the
  weight of $F.$ The argument depends heavily on the decomposition of $V_{\Z\alpha}^+$ as an $M(1)^+$-module and the fusion rules.

We now explain why we only consider  the characterization of
$V_{\Z\alpha}^+$ with $(\alpha,\alpha)/2$ not being a perfect square
in the present paper. In this case, $V_{\Z\alpha}^+$ contains a
primary vector whose weight is not a perfect square. On the other
hand, if $(\alpha,\alpha)/2$ is a perfect square, the weight of any
primary vector in $V_{\Z\alpha}^+$ is a
 perfect square. One assumption in the paper is that there is a primary vector whose weight is not a perfect square.
This assumption is crucial in producing a vertex operator subalgebra
isomorphic to $M(1)^+$. Without this assumption, it is much more
difficult to obtain a subalgebra inside $V$ isomorphic to $M(1)^+$
as the fusion rules among irreducible $L(1,0)$-modules $L(1,n^2)$
are more complicated (see Theorem \ref{co2.2}).

As in \cite{DJ} there is an assumption that $V$ is a sum of highest
weight modules for the Virasoro algebra in the paper. Although we
believe that this assumption is not necessary but we do not know how
to obtain this assumption from the others. On the other hand, we do
not need to assume that the effective central charge $\tilde{c}$ is
equal to the central charge $c$ in this paper as $V_2$ is assumed to
be one dimensional.

This paper is organized as follows. In Section 2 we recall some
basic and important results on vertex operator algebras $L(1,0),$
$M(1)^+$ and $V_L^+$ including the fusion rules. Section 3 is
 about the ${\cal W}_3$ algebra. We show that the simple vertex operator algebra associated to ${\cal W}_3$ of
  central charge $1$ is not a completely reducible module for the Virasoro algebra. This result will be
  used in later sections to deal with the case that the dimension of the weight 3 subspace of the vertex
  operator algebra is greater than 2. In Section 4 we demonstrate that the
  vertex operator algebra $V$ contains a subalgebra isomorphic to $M(1)^+$. In Section 5, we first prove that  $V$  is a completely reducible
   $M(1)^+$-module. Then we give the main result of this paper. That is, a vertex operator algebra satisfying
    certain conditions is isomorphic to $V_L^+$ for some rank one lattice $L.$

\section{Preliminaries}
\def\theequation{2.\arabic{equation}}
\setcounter{equation}{0}

 Let $V=(V,Y,{\bf 1},\omega)$ be a vertex operator algebra \cite{B},
\cite{FLM}. We assume the authors are familiar with  various notions
of $V$-modules  and the definition of rational vertex operator
algebras (cf. \cite{FLM}, \cite{Z}, \cite{DLM}).  We briefly review
the vertex operator algebras associated to the highest weight
representations for the Virasoro algebra and the rank one rational
vertex operator algebras $V_{L}^{+}$.

Here are some basic facts about the highest weight modules for the
Virasoro algebra $Vir$. Let $c,h\in\C$ and $V(c,h)$ be the highest
weight module for the Virasoro algebra $Vir$ with central charge $c$
and highest weight $h.$ We set $\bar{V}(c,0)=V(c,0)/U(Vir)L_{-1}v$
where $v$ is a highest weight vector with highest weight 0 and
denote the irreducible quotient of $V(c,h)$ by $L(c,h).$ We have
(see \cite{KR}, \cite{FZ}):

\begin{prop}\label{vir} Let $c$ be a complex number.

(1) $\bar{V}(c,0)$ is a vertex operator algebra and $L(c,0)$ is a
simple vertex operator algebra.

(2) For any $h\in{\mathbb C}$,  $V(c,h)$ is a module for $\bar
V(c,0).$

(3) $V(c,h)=L(c,h)$, ${\bar V}(c,0)=L(c,0)$,  for $c>1$ and $h>0.$

(4) $V(1,h)=L(1,h)$ if and only if $h\ne \frac{m^2}{4}$ for all
$m\in \Z.$ In case $h=m^2$ for a nonnegative integer $m,$ the unique
maximal submodule of $V(1,m^2)$ is generated by a highest weight
vector with highest weight $(m+1)^2$ and  is isomorphic to
$V(1,(m+1)^2).$
\end{prop}

We need to review the vertex operator algebras $M(1)^+,$ $V_L^+$
and related results \cite{A1}, \cite{A2}, \cite{AD}, \cite{ADL},
\cite{DN1}, \cite{DN2}, \cite{DN3}, \cite{DJL},
 \cite{FLM}.

Let $L=\Z \alpha$ be a positive definite even lattice of rank one.
That is, $(\alpha,\alpha)=2k$ for some positive integer $k.$  Set
$\h=\C\otimes_{\Z} L$ and extend $(\cdot\,,\cdot)$ to a
$\C$-bilinear form on $\h$. Let
$\hat{\h}=\C[t,t^{-1}]\otimes\h\oplus\C K$ be the affine Lie algebra
associated to the abelian Lie algebra $\h$ so that
\begin{align*}
[\alpha(m),\,\alpha(n)]=2km\delta_{m,-n}K\hbox{ and }[K,\hat{\h}]=0
\end{align*}
for any $m,\,n\in\Z$ where $\alpha(m)=\alpha\otimes t^m.$  Then
$\hat{\h}^{\geq 0}=\C[t]\otimes\h\oplus\C K$ is a commutative
subalgebra. For any $\lambda\in\h$, we  define a one dimensional
$\hat{\h}^{\geq 0}$-module $\C e^\lambda$ such that $\alpha(m)\cdot
e^\lambda=(\lambda,\alpha)\delta_{m,0}e^\lambda$ and $K\cdot
e^\lambda=e^\lambda$ for $m\geq0$. We denote by
\begin{align*}
M(1,{\lambda})=U(\hat{\h})\otimes_{U(\hat{\h}^{\geq 0})}\C
e^\lambda\cong S(t^{-1}\C[t^{-1}])
\end{align*}
the $\hat{\h}$-module induced from $\hat{\h}^{\geq 0}$-module $\C
e^\lambda$. Set
$$M(1)=M(1,0).$$ Then there exists a linear map
$Y:M(1)\to\End M(1)[[z,z^{-1}]]$ such that $(M(1),\,Y,\,\1,\,\w)$
carries a simple vertex operator algebra structure and
$M(1,\lambda)$ becomes an irreducible $M(1)$-module for any
$\lambda\in\h$ (see \cite{FLM}). The vacuum vector and the Virasoro
element are given by $\1=e^0$ and $\w=\frac{1}{4k}\alpha(-1)^2\1,$
respectively.

Let $\C[L]$ be the group algebra of $L$ with a basis $e^{\beta}$ for
$\beta\in L.$ The lattice vertex operator algebra associated to $L$
is given by
$$V_L=M(1)\otimes \C[L].$$
The dual lattice $L^{\circ}$ of $L$ is
$$L^{\circ}=\{\,\lambda\in\h\,|\,(\alpha,\lambda)\in\Z\,\}=\frac{1}{2k}L.$$
Then $L^{\circ}=\cup_{i=-k+1}^k(L+\lambda_i)$ is the coset
decomposition with $\lambda_i=\frac{i}{2k}\alpha.$ In particular,
$\lambda_0=0.$ Set $\C[L+\lambda_i]=\bigoplus_{\beta\in L}\C
e^{\beta+\lambda_i}.$ Then each $\C[L+\lambda_i]$ is an
$L$-submodule in an obvious way. Set
$V_{L+\lambda_i}=M(1)\otimes\C[L+\lambda_i]$. Then $V_L$ is a
rational vertex operator algebra and $V_{L+\lambda_i}$ for
$i=-k+1,...,k$ are the irreducible modules for $V_L$ (see \cite{B},
\cite{FLM}, \cite{D1}, \cite{DLM}).

Define a linear isomorphism $\theta:V_{L+\lambda_i}\to
V_{L-\lambda_i}$ for $i\in\{-k+1,...,k\}$ by
\begin{align*}
\theta(\alpha(-n_{1})\alpha(-n_{2})\cdots \alpha(-n_{k})\otimes
e^{\beta+\lambda_i})=(-1)^{k}\alpha(-n_{1})\alpha(-n_{2})\cdots
\alpha(-n_{k})\otimes e^{-\beta-\lambda_i}
\end{align*}
where $n_j>0$ and $\beta\in L.$  Then $\theta$ defines a linear
isomorphism from $V_{L^{\circ}}=M(1)\otimes \C[L^{\circ}]$ to itself
such that
 $$\theta(Y(u,z)v)=Y(\theta u,z)\theta v$$
for $u\in V_{L}$ and $v\in V_{L^{\circ}}.$ In particular, $\theta$
is an automorphism of $V_{L}$ which induces an automorphism of
$M(1).$

For any $\theta$-stable subspace $U$ of $V_{L^{\circ}}$, let $U^\pm$
be the $\pm1$-eigenspace of $U$ for $\theta$. Then $V_L^+$ is a
simple vertex operator algebra.

Also recall the $\theta$-twisted Heisenberg algebra $\h[-1]$ and its
irreducible module $M(1)(\theta)$ from \cite{FLM}. Let $\chi_s$ be a
character of $L/2L$ such that $\chi_s(\alpha)=(-1)^s$ for $s=0,1$
and $T_{\chi_s}=\C$ the irreducible $L/2L$-module with character
$\chi_s$. It is well known that
$V_L^{T_{\chi_s}}=M(1)(\theta)\otimes T_{\chi_s}$ is an irreducible
$\theta$-twisted $V_L$-module (see \cite{FLM}, \cite{D2}). We define
actions of $\theta$ on  $M(1)(\theta)$ and $V_L^{T_{\chi_s}}$ by
\begin{align*}
\theta(\alpha(-n_{1})\alpha(-n_{2})\cdots
\alpha(-n_{k}))=(-1)^{k}\alpha(-n_{1})\alpha(-n_{2})\cdots
\alpha(-n_{k})
\end{align*}
\begin{align*}
\theta(\alpha(-n_{1})\alpha(-n_{2})\cdots \alpha(-n_{k})\otimes
t)=(-1)^{k}\alpha(-n_{1})\alpha(-n_{2})\cdots \alpha(-n_{k})\otimes
t
\end{align*}
for $n_j\in \frac{1}{2}+\Z_{+}$ and $t\in T_{\chi_s}$. We denote the
$\pm 1$-eigenspaces of $M(1)(\theta)$ and $V_L^{T_{\chi_s}}$ under
$\theta$ by $M(1)(\theta)^{\pm}$ and $(V_L^{T_{\chi_s}})^{\pm}$
respectively. We have the following results:
\begin{theorem}\label{t32}
Any irreducible module for the vertex operator algebra $M(1)^+$ is
isomorphic to one of the following modules:$$ M(1)^+, M(1)^-, M(1,
\lambda) \cong M(1, -\lambda)\ (0\neq \lambda \in \h),
M(1)(\theta)^+, M(1)(\theta)^- .$$
\end{theorem}
\begin{theorem}\label{t33}
Any irreducible $V_L^+$-module is isomorphic to one of the following
modules:
$$V_L^{\pm}, V_{\lambda_i+L}( i \not= k),
V_{\lambda_k+L}^{\pm}, (V_L^{T_{\chi_s}})^{\pm}.$$
\end{theorem}

\begin{theorem}\label{t34}  $V_{L}^{+}$ is rational.
\end{theorem}

We remark that the classification of irreducible modules for
arbitrary $M(1)^+$ and $V_L^+$ are obtained in \cite{DN1}-\cite{DN3}
and \cite{AD}. The rationality of $V_L^+$ is established in
\cite{A2} for rank one lattice  and \cite{DJL} in general.

\vskip 0.3cm We next turn our attention to the fusion rules of
vertex operator algebras. Let $V$ be a vertex operator algebra, and
$ W^i$ $ (i=1,2,3$)
 be  ordinary $V$-modules. We denote by $I_{V} \left(\hspace{-3 pt}\begin{array}{c} W^3\\
W^1\,W^2\end{array}\hspace{-3 pt}\right)$  the vector space of all
intertwining operators of type $\left(\hspace{-3 pt}\begin{array}{c}
W^3\\ W^1\,W^2\end{array}\hspace{-3 pt}\right)$.
 For a $V$-module $W$, let
$W^{\prime}$ denote the graded dual of $W$. Then $W'$ is also a
$V$-module \cite{FHL}. It is well known that fusion rules have the
following symmetry (see \cite{FHL}).

\begin{prop}\label{p4.2}
Let $W^{i}$ $(i=1,2,3)$ be $V$-modules. Then
$$\dim I_{{V}} \left(\hspace{-3 pt}\begin{array}{c} W^3\\
W^1\,W^2\end{array}\hspace{-3 pt}\right)=\dim I_{{V}} \left(\hspace{-3 pt}\begin{array}{c} W^3\\
W^2\,W^1\end{array}\hspace{-3 pt}\right), \ \ \ \dim I_{{V}} \left(\hspace{-3 pt}\begin{array}{c} W^3\\
W^1\,W^2\end{array}\hspace{-3 pt}\right)=\dim I_{{V}} \left(\hspace{-3 pt}\begin{array}{c} (W^2)^{\prime}\\
W^1\,(W^3)^{\prime}\end{array}\hspace{-3 pt}\right).$$
\end{prop}

The following two results were obtained in \cite{M} and \cite{DJ}.
\begin{theorem}\label{co2.2} (1) We have
$$
\dim I_{L(1,0)} \left(\hspace{-3 pt}\begin{array}{c} L(1,k^{2})\\
L(1, m^{2})\,L(1, n^{2})\end{array}\hspace{-3 pt}\right)=1,\ \
k\in{\mathbb Z}_{+},  \ |n-m|\leq k\leq n+m,$$
$$\dim I_{L(1,0)} \left(\hspace{-3 pt}\begin{array}{c} L(1,k^{2})\\
L(1, m^{2})\,L(1, n^{2})\end{array}\hspace{-3 pt}\right)=0,\ \
k\in{\mathbb Z}_{+},  \ k<|n-m| \ {\rm or} \  k>n+m, $$ where
$n,m\in{\mathbb Z}_{+}$.

(2) For $n\in{\mathbb Z}_{+}$ such that $n\neq p^{2}$, for all
$p\in{\mathbb Z}_{+}$, we have
$$\dim I_{L(1,0)} \left(\hspace{-3 pt}\begin{array}{c} L(1,n)\\
L(1, m^{2})\,L(1, n)\end{array}\hspace{-3 pt}\right)=1,$$
$$\dim I_{L(1,0)} \left(\hspace{-3 pt}\begin{array}{c} L(1,k)\\
L(1, m^{2})\,L(1, n)\end{array}\hspace{-3 pt}\right)=0,$$
 for $k\in{\mathbb Z}_{+}$ such that  $k\neq n$.

(3)  Let $U$ be a highest weight module for the Virasoro algebra
generated by the highest weight vector $u^{(r)}$ such that
$$
L(0)u^{(r)}=r^{2}u^{(r)}, \ L(k)u^{(r)}=0, \ k\in{\mathbb
Z}_{+}\setminus\{0\}.$$ Let $m,n\in{\mathbb Z}_{+}\setminus \{0\}$
be such that $m\neq n$ and $m,n$ are not perfect squares. Then
$$
I_{L(1,0)} \left(\hspace{-3 pt}\begin{array}{c} U\\
L(1, m)\,L(1, n)\end{array}\hspace{-3 pt}\right)=0.
$$
\end{theorem}

We also need the following result  from \cite{A1} later on.
\begin{theorem}\label{tt33}
Let $M$,$N$ and $T$ be irreducible $M(1)^{+}$-modules. If
$M=M(1,\la)$ such that $\la\neq 0$, then
 $$\dim I_{M(1)^{+}} \left(\hspace{-3 pt}\begin{array}{c} T\\
M\,N\end{array}\hspace{-3 pt}\right)=0 \ \ {\rm or} \ \ 1$$ and
$$\dim I_{M(1)^{+}} \left(\hspace{-3 pt}\begin{array}{c} T\\
M\,N\end{array}\hspace{-3 pt}\right)=1$$ if and only if $(N,T)$ is
one of the following pairs:
$$
(M(1)^{\pm},M(1,\mu))(\la^{2}=\mu^{2}), \ (M(1,\mu),M(1,\nu)), \
(\nu^{2}=(\la\pm \mu)^{2}),$$$$
(M(1)(\theta)^{\pm},M(1)(\theta)^{\pm}), \
(M(1)(\theta)^{\pm},M(1)(\theta)^{\mp}).
$$
\end{theorem}

\section{${\mathcal W}_{3}$ algebra}
\def\theequation{3.\arabic{equation}}
\setcounter{equation}{0} In this section, we recall  ${\mathcal
W}_{3}$ algebra, the associated vertex operator algebra
 and its irreducible quotient (see
\cite{BMP}, \cite{W} and references therein). We will also give some
new results on ${\mathcal W}_3$ with the central charge $C=1$.

 Let ${\mathcal W}_{3}$ be the
associative algebra generated by $L_{m}, W_{m}, m\in \Z$ subject to
the following relations:
$$
[L_{m},L_{n}]=(m-n)L_{m+n}+\delta_{m+n,0}\frac{m^{3}-m}{12}C,
$$
\begin{equation}\label{3e1}
[L_{m},W_{n}]=(2m-n)W_{m+n},
\end{equation}
$$
[W_{m},W_{n}]=(m-n)[\frac{1}{15}(m+n+2)(m+n+3)-\frac{1}{6}(m+2)(n+2)]L_{m+n}$$
\begin{equation}\label{3e2}
+\frac{16}{22+5C}(m-n)\Lambda_{m+n}+\frac{m(m^2-1)(m^2-4)}{360}\delta_{m+n,0}C,\end{equation}
where $$\Lambda_{m}=\sum\limits_{k\geq
2}L_{-k}L_{m+k}+\sum\limits_{k\geq
-1}L_{m-k}L_{k}-\frac{3}{10}(m+2)(m+3)L_{m},$$ and $C$ is a non-zero
central element.  Denote by ${\mathcal W}_{3,\pm}$ and ${\mathcal
W}_{3,0}$ the subalgebras of ${\mathcal W}_{3}$ generated by
$\{L_{m},W_{m} | \pm m>0\}$ and $\{L_{0},W_{0}, C\},$ respectively.
Then ${\mathcal W}_3^{\geq 0}={\mathcal W}_{3,+}+{\mathcal W}_{3,0}$
is also a subalgebra. For $c, \la,\mu\in\C$, let
$\C_{\la,\mu}=\C{\bf 1}_{\la,\mu}$ be the one dimensional module of
${\mathcal W}_{3}^{\geq0}$ such that
$$
L_{0}{\bf 1}_{\la,\mu}=\la{\bf 1}_{\la,\mu}, \ W_{0}{\bf
1}_{\la,\mu}=\mu{\bf 1}_{\la,\mu}, $$
$$ C\cdot {\bf
1}_{\la,\mu}=c{\bf 1}_{\la,\mu}, \  L_m{\bf 1}_{\la,\mu}=W_m{\bf
1}_{\la,\mu}=0$$ for $m>0.$ Denote by ${\mathcal M}(c,\la,\mu)$ the
induced module:
$${\mathcal M}(c,\la,\mu)={\mathcal
W}_3\otimes_{{\mathcal W}_{3}^{\geq 0}}\C_{\la,\mu}.$$

It is well known that ${\mathcal M}(c,\la,\mu)$ has a unique
irreducible quotient which is denoted by ${\mathcal L}(c,\la,\mu)$.
It is easy to see that
$$\C L_{-1}{\bf 1}+\C W_{-1}{\bf 1}+\C W_{-2}{\bf 1}$$
is an invariant subspace of ${\mathcal M}(c,0,0)$ under the action
of ${\mathcal W}_{3,+}.$  Let ${\mathcal J}$ be the submodule of
${\mathcal M}(c,0,0)$ generated by the  three vectors $L_{-1}{\bf
1}, W_{-1}{\bf 1}$ and $W_{-2}{\bf 1}$. Denote by $\bar{\mathcal
M}(c,0,0)$ the quotient ${\mathcal M}(c,0,0)/{\mathcal J}$. Let
$$
L(z)=\sum\limits_{n\in\Z}L_{n}z^{-n-2},  \
W(z)=\sum\limits_{n\in\Z}W_{n}z^{-n-3}.
$$
Then $\bar{\mathcal M}(c,0,0)$ is a vertex operator algebra
generated by $\omega=L_{-2}{\bf 1}$ and $w=W_{-3}{\bf 1}$ such that
$Y(\omega, z)=L(z)$ and $Y(w,z)=W(z).$ The irreducible quotient
${\mathcal L}(c,0,0)$ is the associated simple vertex operator
algebra. The following lemma is clear (see \cite{BMP}, \cite{W}).
\begin{lem}\label{3l1}
$\bar{\mathcal M}(c,0,0)$ has a linear basis:
$$
L_{-m_{1}}L_{-m_{2}}\cdots L_{-m_{s}}W_{-n_{1}}W_{-n_{2}}\cdots
W_{-n_{t}}{\bf 1},
$$
where $s,t\geq 0$, $m_{1}\geq m_{2}\geq\cdots \geq  m_{s}\geq 2$,
$n_{1}\geq n_{2}\geq\cdots\geq n_{t}\geq 3.$
\end{lem}

In the following discussion, we assume that the central charge
$c=1$. Note that the vertex operator subalgebra
$\langle\omega\rangle$ is isomorphic to the irreducible module
$L(1,0)$ for the Virasoro algebra. In the following we will identify
$L(1,h)$ with its isomorphic image in $\bar{\mathcal M}(1,0,0)$ if
there is no confusion arising.

For convenience, we call a non-zero element $v$ in $\bar{\mathcal
M}(1,0,0)$ a primary vector if $L_{m}v=0$, for all $m\geq 1$. From
the definition of $\bar{\mathcal M}(1,0,0)$, it is obvious that
there is no primary element of weight $4$ or $5$, and
$$\dim \bar{\mathcal M}(1,0,0)_{r}=r-1,
\ 3\leq r\leq 5.$$ By Lemma \ref{3l1}, $\bar{\mathcal M}(1,0,0)_{6}$
has a basis:
$$
W_{-3}W_{-3}{\bf 1}, \ W_{-6}{\bf 1}, \ L_{-2}W_{-4}{\bf 1}, \
L_{-3}W_{-3}{\bf 1}, \ L_{-6}{\bf 1}, \ L_{-4}L_{-2}{\bf 1}, \
L_{-3}^{2}{\bf 1}, \ L_{-2}^{3}{\bf 1}.$$ The $L(1,0)$-submodule of
$\bar{\mathcal M}(1,0,0)$ generated by $w=W_{-3}{\bf 1}$ is
isomorphic to $L(1,3).$ Note that
$$
\dim L(1,0)_{6}=4, \ \ \dim (L(1,3)\cap \bar{\mathcal
M}(1,0,0)_{6})=3,$$ then there is a non-zero element of weight $6$
in $\bar{\mathcal M}(1,0,0)$ which is not in $L(1,0)\oplus L(1,3)$.
Since $\dim (\bar{\mathcal M}(1,0,0)_{4}+\bar{\mathcal
M}(1,0,0)_{5})=7<8$, it follows that  there is a primary element
$u^{(6)}$ of weight $6$ in $\bar{\mathcal M}(1,0,0)$. Assume that
$$
u^{(6)}=a_{1}W_{-3}w+
a_{2}L_{-3}w+a_{3}L_{-2}L_{-1}w+a_{4}L_{-1}^{3}w+a_{5}v,$$ for some
$v\in L(1,0)$ and $a_{i}\in \C$, $1\leq i\leq 5$. From the
commutation relation (\ref{3e2}) we see that
$$L_{i}W_{-3}w\in L(1,0), \ i\geq 1.$$
Note that
$$
L_{i}(a_{2}L_{-3}w+a_{3}L_{-2}L_{-1}w+a_{4}L_{-1}^{3}w)\in L(1,3)$$
for $i\geq 1.$ This implies that $a_{i}=0$ for $i=2,3,4$.
Consequently, $a_{1}\neq 0$ and
\begin{equation}\label{2e}
W_{-3}w\in L(1,0)\oplus L(1,6).\end{equation}

The $L(1,0)$-submodule  generated by $u^{(6)}$ is isomorphic to
$L(1,6).$ By computing the dimensions of $L(1,0)_7, L(1,3)_7,
L(1,6)_7$ and $\bar{\mathcal M}(1,0,0)_{7}$ we immediately see that
there is no primary vector in $\bar{\mathcal M}(1,0,0)_{7}$. Note
that $\bar{\mathcal M}(1,0,0)_{8}$ has a basis:
$$
\{W_{-5}w,\ W_{-4}W_{-4}{\bf 1}, \ L_{-2}W_{-3}w,  \ W_{-8}{\bf 1},
\ L_{-2}W_{-6}{\bf 1}, \ L_{-3}W_{-5}{\bf 1}, $$$$\ L_{-4}W_{-4}{\bf
1}, \ L_{-2}^{2}W_{-4}{\bf 1}, L_{-5}w, L_{-3}L_{-2}w\}\cup S_{8},$$
and the weight 8 subspace of $L(1,6)\oplus L(1,3)\oplus L(1,0)$ has
a basis:
$$
\{L_{-2}W_{-3}w, \ L_{-1}^{2}W_{-3}w, \ L_{-5}w, \ L_{-4}L_{-1}w, \
L_{-3}L_{-2}w, $$$$ \ L_{-3}L_{-1}^{2}w, \ L_{-2}^{2}L_{-1}w, \
L_{-2}L_{-1}^{3}w, \ L_{-1}^{5}w\}\cup S_{8},
$$ where $S_{8}$ is a basis of $L(1,0)_{8}$. It follows that there exists a
non-zero element  in $\bar{\mathcal M}(1,0,0)_{8}$ which is not in
the $L(1,0)$-submodule generated by $u^{(6)}$, $w$ and ${\bf 1}$.

We are now in a position to state the main result of this section.
\begin{theorem}\label{3t1}
The simple vertex operator algebra ${\mathcal L}(1,0,0)$ is not
completely reducible as an $L(1,0)$-module.
\end{theorem}
\pf By the above discussion, we have
$$\dim \bar{\mathcal M}(1,0,0)_{8}-\dim (L(1,6)\oplus
L(1,3)\oplus L(1,0))_{8}=1. $$ Then there exists
$u^{(8)}\in\bar{\mathcal M}(1,0,0)_{8}$ such that $u^{(8)}$ is
either a primary element or the irreducible quotient of
$L(1,0)$-module generated by $u^{(8)}, u^{(6)},w$ and $\omega$
modulo the submodule $L(1,6)\oplus L(1,3)\oplus L(1,0)$ is
isomorphic to $L(1,8)$. It is easy to check that $u^{(8)}$,
$u^{(6)}$ and $w$  do not lie in the maximal submodule $J$.

Now let
\begin{align*}
u^{(9)} =& 3312738 W_{-3}W_{-3}w-214776 L_{-6}w-4379460 L_{-5}L_{-1}w-10304064  L_{-4}L_{-2}w\\
&
-7494075 L_{-3}^{2}w+2682708 L_{-4}L_{-1}^{2}w+8127252 L_{-3}L_{-2}L_{-1}w+424032 L_{-2}^{3}w\\
& -2068431 L_{-3}L_{-1}^{3}w-2594664 L_{-2}^{2}L_{-1}^{2}w+350889
L_{-2}L_{-1}^{4}w+159578 L_{-1}^{6}w
\end{align*}
A direct calculation yields that
$$
L(n)u^{(9)}=0,\ \ n\geq 1.$$ Furthermore, we have
\begin{align*}
 W_{1}W_{-3}W_{-3}w=&\frac{146}{3}L_{-2}W_{-3}w+\frac{10070}{27}L_{-8}{\bf
 1}+\frac{49664}{729}L_{-5}L_{-3}{\bf 1}
 +\frac{20480}{729}L_{-3}L_{-2}{\bf
 1}\\
 &+\frac{7232}{729}L_{-4}^{2}{\bf
 1}+\frac{28672}{729}L_{-4}L_{-2}^{2}{\bf 1}
 +\frac{95872}{729}L_{-6}L_{-2}{\bf 1}+\frac{310}{9}W_{-5}w,
 \end{align*}
$$W_{1}L_{-6}w=13W_{-5}w+2L_{-6}L_{-2}{\bf 1},\\
W_{1}L_{-5}L_{-1}w=11W_{-4}W_{-4}{\bf 1}+5L_{-5}L_{-3}{\bf 1},$$
$$W_{1}L_{-4}L_{-2}w=9L_{-2}W_{-3}w+9W_{-5}w-\frac{5}{9}L_{-4}^{2}{\bf
1}+\frac{214}{27}L_{-4}L_{-2}^{2}{\bf 1},$$
$$W_{1}L_{-3}^{2}w=28W_{-5}w-\frac{56}{27}L_{-8}{\bf
1}-\frac{28}{27}L_{-5}L_{-3}{\bf
1}+\frac{502}{27}L_{-3}^{2}L_{-2}{\bf 1},$$
\begin{align*}
W_{1}L_{-4}L_{-1}^{2}w=&-10L_{-8}{\bf
1}+18W_{-5}w+\frac{128}{3}L_{-6}L_{-2}{\bf
1}+\frac{128}{3}L_{-5}L_{_3}{\bf 1}\\
&+\frac{110}{3}L_{-4}^{2}{\bf 1}+\frac{64}{9}L_{-4}L_{-2}^{2}{\bf
1},\\
W_{1}L_{-3}L_{-2}L_{-1}w=&14W_{-4}W_{-4}{\bf 1}+\frac{704}{27}L_{-8}{\bf 1}+\frac{280}{9}L_{-6}L_{-2}{\bf 1}\\
& +\frac{448}{27}L_{-4}L_{-2}^{2}{\bf
1}+\frac{688}{27}L_{-5}L_{-3}{\bf
1}+\frac{839}{27}L_{-3}^{2}L_{-2}{\bf 1},\\
W_{1}L_{-2}^{3}w=& 45L_{-2}W_{-3}w+15W_{-5}w-\frac{40}{3}L_{-8}{\bf
1}\\
&-\frac{20}{3}L_{-6}L_{-2}{\bf 1}-\frac{5}{3}L_{-4}L_{-2}^{2}{\bf
1}+\frac{178}{9}L_{-2}^{4}{\bf 1},\\
W_{1}L_{-3}L_{-1}^{3}w=&\frac{1064}{9}L_{-8}{\bf
1}+\frac{1792}{9}L_{-6}L_{-2}{\bf 1}+\frac{2520}{9}L_{-5}L_{-3}{\bf
1}\\ &+\frac{896}{9}L_{-4}^{2}{\bf
1}+\frac{448}{9}L_{-3}^{2}L_{-2}{\bf
1},\\
W_{1}L_{-2}^{2}L_{-1}^{2}w=&30W_{-5}w+\frac{8462}{27}L_{-8}{\bf
1}+\frac{3152}{9}L_{-6}L_{-2}{\bf 1}+\frac{4480}{27}L_{-5}L_{-3}{\bf
1}+\frac{320}{9}L_{-4}^{2}{\bf 1}\\
&+\frac{2974}{27}L_{-4}L_{-2}^{2}{\bf
1}+\frac{1280}{27}L_{-3}^{2}L_{-2}{\bf 1}+\frac{64}{9}L_{-2}^{4}{\bf
1},\\
W_{1}L_{-2}L_{-1}^{4}w=&\frac{25144}{9}L_{-8}{\bf
1}+\frac{5280}{3}L_{-6}L_{-2}{\bf 1}+\frac{9728}{9}L_{-5}L_{-3}{\bf
1}\\
&+\frac{1280}{3}L_{-4}^{2}{\bf 1}+\frac{2048}{9}L_{-4}L_{-2}^{2}{\bf
1}+\frac{1024}{9}L_{-3}^{2}L_{-2}{\bf 1},\\
W_{1}L_{-1}^{6}w=&26800L_{-8}{\bf 1}+\frac{25600}{3}L_{-6}L_{-2}{\bf
1}+\frac{25600}{3}L_{-5}L_{-3}{\bf 1}+\frac{12800}{3}L_{-4}^{2}{\bf
1}.
\end{align*}
Note that the monomial $L_{-2}^{4}{\bf 1}$ appears only in
$W_{1}L_{-2}^{3}w$ and $W_{1}L_{-2}^{2}L_{-1}^{2}w$. Then it is easy
to check  that
$$
W_{1}u^{(9)}\neq 0.$$ So $u^{(9)}$ does not lie in $J.$ This implies
that $u^{(9)}$ is a primary element of ${\mathcal L}(1,0,0)$. If
${\mathcal L}(1,0,0)$ is completely reducible as an $L(1,0)$-module,
then it follows from (\ref{2e}) and the expression of $u^{(9)}$
that
  the space of intertwining operators
of type
$$\left(\hspace{-3 pt}\begin{array}{c}
L(1,9)\\
L(1,3)\, \ \  L(1,6)\end{array}\hspace{-3 pt}\right)$$ is non-zero,
which is a contradiction with Theorem \ref{co2.2}. Thus ${\mathcal
L}(1,0,0)$ is not completely reducible as an $L(1,0)$-module. \qed

\section{Existence of $M(1)^{+}$ in $V$}
\def\theequation{4.\arabic{equation}}
\setcounter{equation}{0}

From now on, we always assume that $V$ is a simple rational vertex
operator algebra of central charge 1 satisfying the following
conditions:

(1) \ $V=\oplus_{n=0}^{\infty}V_{n}, \ V_{0}=\C {\bf 1}, \ V_{1}=0,
\ \dim V_{2}=1;$

(2) \ $\dim V_{3}\geq 2$, or $\dim V_{3}=1$ and $\dim V_{4}\geq 3$;

(3) \ $V$ is a sum of highest weight modules of $L(1,0)$.

\begin{remark}\label{4r1} Let $V$ be a rational,
$C_{2}$-cofinite  CFT type vertex operator algebras with
$c=\tilde{c}=1.$ It is shown in \cite{DM1} that if $V_1\ne 0$ then
$V$ is isomorphic to $V_L$ for a rank one lattice $L.$ It is
established in \cite{ZD} and \cite{DJ} that if $V_1=0$ and $\dim
V_2>1$ then $V$ is isomorphic to $V_{\Z\alpha}^+$ with
$(\alpha,\alpha)=4.$ So for the purpose of the classification we
need only to consider the case that $V_{1}=0$ and $\dim V_{2}=1.$
For the characterization of the rational vertex operator algebra
$V_{\Z\al}^{+}$, it is very natural to have the assumption (2). We
expect the assumption (3) is true for any rational vertex operator
algebra with $c=\tilde{c}=1$, but it seems difficult to prove it
right now.
\end{remark}
\vskip 0.3cm We have the following lemma from \cite{DJ}.
\begin{lem}\label{l1} $V$ is a completely reducible
module for the Virasoro vertex operator  algebra $L(1,0)$.
\end{lem}

It is obvious that $V$ carries a non-degenerate symmetric bilinear
form $(\cdot,\cdot)$ such that $({\bf 1},{\bf 1})=1$ (\cite{FHL},
\cite{L1}).
 We will prove in this
section that $V$ contains a vertex operator subalgebra isomorphic to
$M(1)^{+}$. Let $X^{1}$ and $X^{2}$ be two subsets of $V$. Set
$$
X^{1}\cdot X^{2}={\rm span}\{x_{n}y|x\in X^{1}, y\in X^{2}, n\in{\mathbb
Z}\}.$$
 Recall from Section 3 that a non-zero element $v$ in $V$ is called a primary vector if $L_{n}v=0$ for all $n\geq 1$. We first have the following lemma.
\begin{lem}\label{ll2.4}
 Let $u^{1},u^{2}\in V$ be two primary elements. Let  $U^{1}$ and $U^{2}$ be the  two $L(1,0)$-submodules of $V$ generated
by $u^{1}$ and $u^{2}$  respectively. Then
$$
U^{1}\cdot U^{2}={\rm span}\{L(-m_{1})\cdots L(-m_{s})u^{1}_{n}u^{2} |
s\geq 0, m_{1},\cdots, m_{s}\in {\mathbb Z}_{+}, n\in{\mathbb
Z}\}.$$
\end{lem}
\pf Let $x=L(-m_{1})\cdots L(-m_{s})u^{1}$ and
$u=x_{n}L(-n_{1})\cdots L(-n_{k})u^{2}$. Using the formula
$$(a_lb)_m=\sum_{j\geq 0}(-1)^j{l\choose j}a_{l-j}b_{m+j}
-\sum_{j\geq 0}(-1)^{l+j}{l\choose j}b_{m+l-j}a_{j}$$ (which holds
for any vertex operator algebra and any vectors $a,b$ in the vertex
operator algebra) and the commutator formula
$$[L(p), u^i_q]=(-p-q-1+(p+1)\wt u^i)u^i_{p+q}$$
for $l,p,q\in\Z$ and $i=1,2$  we see that
\begin{eqnarray*}
& & u=x_n(L(-n_{1})\cdots L(-n_{k})u^{2})\\
& &\ \ =\sum_{p_1,...,p_s, p\in \Z}a_{p_1,...,p_s,p}L(p_1)\cdots
L(p_s)u^1_pL(-n_{1})\cdots L(-n_{k})u^{2}
\end{eqnarray*}
lies in
$${\rm span}\{L(-m_{1})\cdots L(-m_{s})u^{1}_{n}u^{2} |
s\geq 0, m_{1},\cdots, m_{s}\in {\mathbb Z}_{+}, n\in{\mathbb
Z}\},$$ where $a_{p_1,...,p_s,p}\in \C.$  The proof is complete.
\qed

\begin{lem}\label{3l3}
$V$ contains a primary element of weight 4.
\end{lem}
\pf If $\dim V_{3}=1$, then the lemma follows from the assumption
(2) and Lemma \ref{l1}  as $\dim L(1,0)_4=2.$ If $\dim V_{3}\geq 2$,
then there exists at least one primary vector of weight 3 in $V$.
Let $F$ be a primary element of weight 3 such that $(F,F)=2$. Set
$$
u=F_{1}F+\frac{2}{3}L(-4){\bf 1}-\frac{64}{9}L(-2)^{2}{\bf 1}.
$$
Then it is easy to check that $L_{n}u=0$, for all $n\geq 1$.

If $u\neq 0$, then the lemma holds. If $u=0$, denote
$W_{n}=\frac{1}{\sqrt{6}}F_{n+2}$, for $n\in\Z$. Then  it is easy to
check that $L_{m}, W_{n}$ for $m,n\in\Z$ satisfy relations
(\ref{3e1}) and (\ref{3e2}) with $C=1$. Then the vertex operator
subalgebra of $V$ generated by $F$ and $\omega$ is a quotient of
$\bar{\mathcal M}(1,0,0)$.
 By Theorem \ref{3t1}, the vertex operator
subalgebra of $V$ generated by $F$ and $\omega$ is not a completely
reducible $L(1,0)$-module, which contradicts Lemma \ref{l1}. \qed

 Now let $J$ be
a primary element of weight 4. We may assume that
\begin{equation}\label{eq4.1}
 (J,J)=54, \ \
\end{equation}
So $J_{7}J=54{\bf 1}.$

Since $V$ is rational, it follows that $V$ contains infinitely many
primary vectors. In this paper we will deal with the case that $V$
contains a primary vector whose weight is not a perfect square.

Let $k$ be the smallest positive integer such that $k$ is not a
perfect square and $V$ contains a primary vector  of weight $k$.
Then $k\geq 3$. Let $U$ be the subalgebra of $V$ generated by
$\omega$ and $J$. We will prove that $U$ is isomorphic to
$M(1)^{+}$. Let $V^{(4)}$ be the irreducible   $L(1,0)$-submodule of
$U$ generated by $J$. Then $V^{(4)}\cong L(1,4)$. By Theorem
\ref{co2.2}, for $n\in\Z$,
$$
\dim I_{L(1,0)} \left(\hspace{-3 pt}\begin{array}{c} L(1,n)\\
V^{(4)}\, \ L(1, k)\end{array}\hspace{-3 pt}\right)\neq 0
$$
if and only if $n=k$. Let $$A_{k}=\{v\in V_{k}|L_{n}v=0, \ n\geq
1\}.$$ For any $0\ne x\in A_{k},$ the $L(1,0)$-submodule $U(Vir)x$
of $V$ generated by $x$ is isomorphic to $L(1,k)$ and
${\rm span}\{u_nv|u\in V^{(4)}, v\in U(Vir)x, n\in \Z\}$ is an irreducible
$L(1,0)$-module isomorphic to $L(1,k)$ again. In particular,
$J_3x\in A_{k}.$ Let $F$ be an eigenvector of $J$ in $A_{k}.$ Let
$V^{(k)}$ be the irreducible $L(1,0)$-submodule of $V$ generated by
$F$. Then
$$
V^{(4)}\cdot V^{(k)}=V^{(k)}.
$$

Let $(V_{L}^{+}, {Y}(\cdot,z))$ be the rank one rational vertex
operator algebra with $L=\Z\al$ such that $(\al,\al)=2k$. Set
$$J^1= h(-1)^4{\bf 1} -2h(-3)h(-1){\bf 1} + \frac{3}{2}h(-2)^2{\bf 1}, \
E=e^{\al}+e^{-\al}$$ where $h=\frac{\alpha}{\sqrt{2k}}.$  Then
\begin{equation}\label{eq4.2}
(J^1,J^1)=54, \ \ (E,E)=2.
\end{equation}
Moreover, $J^1$ and the conformal element $\omega$ of $V_{L}^{+}$
generate the vertex operator subalgebra $M(1)^{+}$. Let $M^{(4)}$
and $M^{(k)}$ be the irreducible  $L(1,0)$-submodules of $V_L^+$
generated by $J^1$ and $E$ respectively. Then there exist
$L(1,0)$-module isomorphisms $\sigma_{1}$: $M^{(4)}\rightarrow
V^{(4)}$ and  $\sigma_{2}$: $M^{(k)}\rightarrow V^{(k)}$ such that
\begin{equation}\label{eq4.6} \sigma_{1}(J^1)=J, \ \sigma_{2}(E)=F
\end{equation} and
$$
(u^1,v^{1})=(\sigma_{1}(u^1),\sigma_{1}(v^1))$$ for $u^1,v^1\in
M^{(4)}$.

We identify the Virasoro vertex operator algebra $L(1,0)$ in $V$ and
$V_{L}^{+}$.  Let
$${\mathcal I}^0(u^1,z)v^1={\mathcal P^0}\circ
Y(u^1,z)v^1$$
 for $u^1,v^1\in M^{(4)}$ be the intertwining operator
of type
$$\left(\hspace{-3 pt}\begin{array}{c}
L(1,0)\\
M^{(4)}\, \ \  M^{(4)}\end{array}\hspace{-3 pt}\right),$$ and
${\mathcal I}^0(\sigma_{1}(u^1),z)\sigma_{1}(v^1)={\mathcal
Q^0}\circ Y(\sigma_{1}(u^1),z)\sigma_{1}(v^1)$ for $u^1,v^1\in
M^{(4)}$ be the intertwining operator of type
$$\left(\hspace{-3 pt}\begin{array}{c}
L(1,0)\\
V^{(4)}\, \ \  V^{(4)}\end{array}\hspace{-3 pt}\right),$$ where
${\mathcal P^0}$ and ${\mathcal Q^0}$ are the projections of
$V_{L}^{+}$ and $V$ to $L(1,0)$ respectively. By
(\ref{eq4.1})-(\ref{eq4.2}) and Theorem \ref{co2.2}, we have
\begin{equation}\label{eq4.3}
{\mathcal I}^0(u^1,z)v^1={\mathcal
I}^0(\sigma_{1}(u^1),z)\sigma_{1}(v^1),
\end{equation}
for  $u^1,v^1\in M^{(4)}$.  Furthermore, we have the following
lemma.

\begin{lem}\label{ll4.2} Replacing $J$ by $-J$ if necessary, we have
\begin{equation}\label{eq2.4}
\sigma_{2}(Y(u^1,z)v^2)= Y(\sigma_{1}(u^1),z)\sigma_{2}(v^2),
\end{equation}
for $u^{1}\in M^{(4)}, v^2\in M^{(k)}.$
\end{lem}
\pf Since $M^{(4)}\cong V^{(4)}\cong L(1,4)$ and $M^{(k)}\cong
V^{(k)}\cong L(1,k)$, we may identify $M^{(4)}$ with $V^{(4)}$
through $\sigma_{1}$ and $M^{(k)}$ with $V^{(k)}$ through
$\sigma_{2}$. So both $Y(u^1,z)v^2$ and
$Y(\sigma_{1}(u^1),z)\sigma_{2}(v^2)$ for $u^{1}\in M^{(4)}, v^2\in
M^{(k)}$ are intertwining operators of type
$$\left(\hspace{-3 pt}\begin{array}{c}
L(1,k)\\
L(1,4)\, \ \  L(1,k)\end{array}\hspace{-3 pt}\right).$$ By Theorem
\ref{co2.2} and (\ref{eq4.6}), we have
$$\sigma_{2}(Y(u^1,z)v^2)=\varepsilon Y(\sigma_{1}(u^1),z)\sigma_{2}(v^2),
$$
for some $\varepsilon\in \C$. By the Jacobi identity, we have
$$
(J^{1}_{7}J^{1})_{-1}E=\sum\limits_{i=0}^{\infty}(-1)^{i}\left(\begin{array}{c}7\\i\end{array}\right)
(J^{1}_{7-i}J^{1}_{-1+i}+J^{1}_{6-i}J^{1}_{i})E,
$$
$$(J_{7}J)_{-1}F=\sum\limits_{i=0}^{\infty}(-1)^{i}\left(\begin{array}{c}7\\i\end{array}\right)
(J_{7-i}J_{-1+i}+J_{6-i}J_{i})F.$$ By (\ref{eq4.2}), we have
$\sigma_{2}((J^{1}_{7}J^{1})_{-1}E)=54F=(J_{7}J)_{-1}F$. On the
other hand,
\begin{eqnarray*}
&&\sigma_2((J_{7}^1J^1)_{-1}E)=\sigma_2((J_{7-i}^1J_{-1+i}^1+J_{6-i}^1J_{i}^1)E)\\
&&\hspace{1cm}=\varepsilon^2(J_{7-i}J_{-1+i}+J_{6-i}J_{i})F\\
&&\hspace{1cm}=\varepsilon^2(J_{7}J)_{-1}F.
\end{eqnarray*}
This forces $\varepsilon=\pm 1.$ If $\varepsilon=1$ we are done.
Otherwise, we may replace $J$ by $-J$ and have the desired result.
\qed

\begin{lem}\label{ll4.1}
 If there exists a primary element $v$ of weight 4 in $V$
 such that $v_{3}F\in \C F$, then $v\in \C J$.
\end{lem}

\pf  Let $X$ be the irreducible $L(1,0)$-submodule of $V$ generated
by $v$. Then we also have
$$X\cdot V^{(k)}=V^{(k)}$$
and $Y(w^{2} ,z)|_{V^{(k)}}$ for  $w^{2}\in X$ is an intertwining
operator of type $\left(\hspace{-3 pt}\begin{array}{c} V^{(k)}\\
X\,V^{(k)}\end{array}\hspace{-3 pt}\right)$ from the assumption.

Notice that $V^{(4)}$ is isomorphic to $X$ as $L(1,0)$-modules. Let
$\varphi$ be an isomorphism from  $V^{(4)}$ to $X$. Then
${Y}(\varphi(w^{1}),z)|_{V^{(k)}}$ for $w^{1}\in V^{(4)}$ is an
intertwining
operator of type $\left(\hspace{-3 pt}\begin{array}{c} V^{(k)}\\
V^{(4)}\,V^{(k)}\end{array}\hspace{-3 pt}\right)$. By Theorem
\ref{co2.2}, there exists some nonzero $c\in\C$ such that
$${Y}(w^{1}-c\varphi(w^{1}),z)|_{V^{(k)}}=0,$$
for all $w^{1}\in V^{(4)}$. By Proposition 11.9 of \cite{DL}, we see
that $w^{1}-c\varphi(w^{1})=0$ for any $w^1\in V^{(4)}.$ This
implies that $\varphi(w^{1})\in V^{(4)}$ and $V^{(4)}=X.$ The lemma
is proved. \qed

By Theorem \ref{co2.2} and the fact that $V_{1}=0$, we may assume
that
$$J_{3}J=x+y,$$
where $x\in L(1,0)$, $y$ is either zero or a primary element of
weight 4. Since
$$ (J_{3}J)_{3}F\in V^{(k)},$$ it follows that $y_{3}F\in V^{(k)}$.
Then by Lemma \ref{ll4.1}, $y\in V^{(4)}$. Thus
\begin{equation}\label{eeq4.17}
J_{3}J\in L(1,0)\oplus V^{(4)}. \end{equation} Now let ${\mathcal
P}^4$ and ${\mathcal Q}^4$ be the projections of $V_{L}^{+}$  and
$V$ to $M^{(4)}$ and $V^{(4)}$ respectively. Let ${\mathcal
I}^4(u^1,z)v^1={\mathcal P^4}\circ Y(u^1,z)v^1$ for $u^1,v^1\in
M^{(4)}$ be the intertwining operator of type
$$\left(\hspace{-3 pt}\begin{array}{c}
M^{(4)}\\
M^{(4)}\, \ \  M^{(4)}\end{array}\hspace{-3 pt}\right),$$ and
${\mathcal I}^4(\sigma_{1}(u^1),z)\sigma_{1}(v^1)={\mathcal
Q^4}\circ Y(\sigma_{1}(u^1),z)\sigma_{1}(v^1)$ for $u^1,v^1\in
M^{(4)}$ be the intertwining operator of type
$$\left(\hspace{-3 pt}\begin{array}{c}
V^{(4)}\\
V^{(4)}\, \ \  V^{(4)}\end{array}\hspace{-3 pt}\right).$$ We have
the following lemma.

\begin{lem}\label{ll4.3}
 For $u^1,v^1\in M^{(4)}$, we have
$$\sigma_{1}{\mathcal I}^4(u^1,z)v^1={\mathcal
I}^4(\sigma_{1}(u^1),z)\sigma_{1}(v^1).$$
\end{lem}

\pf By Theorem \ref{co2.2}, we have
\begin{equation}\label{eq4.5}\sigma_{1}{\mathcal
I}^4(u^1,z)v^1=c{\mathcal I}^4(\sigma_{1}(u^1),z)\sigma_{1}(v^1),
\end{equation}
for some $c\in\C$. Using Lemma \ref{ll4.2} we see that  for
$n\in\Z$,
$$\sum\limits_{i=0}^{\infty}(-1)^{i}\left(\begin{array}{c}3\\i\end{array}\right)
(J_{3-i}J_{n+i}+J_{3+n-i}J_{i})F=\sigma_{2}(\sum\limits_{i=0}^{\infty}(-1)^{i}\left(\begin{array}{c}3\\i\end{array}\right)
(J^{1}_{3-i}J^{1}_{n+i}+J^{1}_{3+n-i}J^{1}_{i})E).$$ From the Jacobi
identity we know that
$$(J_{3}J)_{n}=\sum\limits_{i=0}^{\infty}(-1)^{i}\left(\begin{array}{c}3\\i\end{array}\right)
(J_{3-i}J_{n+i}+J_{3+n-i}J_{i}).$$ So
$$(J_{3}J)_{n}F=\sigma_{2}((J^1_{3}J^1)_{n}E).$$
That is, $$
((\sigma_{1}J^1)_{3}(\sigma_{1}J^1))_{n}(\sigma_{2}E)=\sigma_{2}((J^1_{3}J^1)_{n}E).$$

From the discussion above we may assume that
 $$J_{3}J=x+y, \ J^1_{3}J^1=x^1+y^1,$$ where
$x,x^1\in L(1,0)$, $y^1\in M^{(4)}$, $y\in V^{(4)}$. Then
$$
(x+y)_{n}F=\sigma_{2}((x^1+y^1)_{n}E)=(\sigma_1(x^1+y^1))_n\sigma_2(E)=(\sigma_1(x^1+y^1))_nF
$$
for all $n.$ Thus $\sigma_1(x^1+y^1)=x+y,$ $x=\sigma_1(x^1)$ and
$y=\sigma_1(y^1).$ On the other hand,  by (\ref{eq4.5}), we have
$\sigma_{1}(y^1)=cy.$ This forces $c=1.$
 \qed

By the skew-symmetry, we have
$$
Y(J,z)J=e^{L(-1)z}Y(J,-z)J.$$It follows that
$$
J_{-2}J=-J_{-2}J+\sum_{j=1}^{9}(-1)^{j+1}\frac{1}{j!}L(-1)^{j}J_{-2+j}J.
$$
This together with Theorem \ref{co2.2} and Lemma \ref{ll2.4} deduces
that
$$
V^{(4)}\cdot V^{(4)}\subseteq L(1,0)\oplus L(1,4)\oplus
L(1,4^{2}).$$ Actually we have

\begin{lem}\label{ll4.4}
 $\ V^{(4)}\cdot V^{(4)}\cong L(1,0)\oplus L(1,4)\oplus
L(1,4^{2}).$
\end{lem}

\pf  Again by Lemma \ref{ll4.2} we see that for $n\in\Z$,
$$\sum\limits_{i=0}^{\infty}(-1)^{i}{ -9\choose i}
(J_{-9-i}J_{n+i}+J_{-9+n-i}J_{i})F=\sigma_{2}(\sum\limits_{i=0}^{\infty}(-1)^{i}{
-9\choose i} (J^{1}_{-9-i}J^{1}_{n+i}+J^{1}_{-9+n-i}J^{1}_{i})E).$$
It follows from the Jacobi identity that
$$(J_{-9}J)_{n}F=\sigma_{2}((J^1_{-9}J^1)_{n}E).$$
That is, $$
((\sigma_{1}J^1)_{-9}(\sigma_{1}J^1))_{n}(\sigma_{2}E)=\sigma_{2}((J^1_{-9}J^1)_{n}E).$$
Let $M^{(4^2)}$ be the irreducible $L(1,0)$-submodule of $M(1)^+$
isomorphic to $L(1,16).$
 Let $$J^1_{-9}J^1=x^0+x^4+x^{4^2}, \ J_{-9}J=y^0+y^4+y^{4^2},$$ where
$x^0,y^0\in L(1,0)$, $x^4\in M^{(4)},$  $y^4\in V^{(4)},$
$x^{4^2}\in M^{(4^2)}$ and $y^{4^2}\in V$ are primary vectors of
weight 16. Then
$$(y^0+y^4+y^{4^2})_{n}F=\sigma_{2}((x^0+x^4+x^{4^2})_{n}E).
$$
By (\ref{eq4.3}), we have
$$x^0=y^0.$$
 By Lemma \ref{ll4.2} and Lemma \ref{ll4.3}, we have  \
\begin{equation}\label{eq4.9} \sigma_{2}(x^4_{n}E)=(\sigma_{1}x^4)_{n}F=
y^4_{n}F. \end{equation}
 So
$$y^{4^2}_{n}F=\sigma_{2}(x^{4^2}_{n}E).$$
Since $M^{(4)}\cdot M^{(4)}\cong L(1,0)\oplus L(1,4)\oplus L(1,16),$
we conclude that $x^{4^2}_{n}E\ne 0.$ This can also be verified
directly. As a result, $y^{4^2}_{n}F\ne 0.$ The lemma  follows. \qed

Denote by $V^{(4^2)}$  the irreducible  $L(1,0)$-submodule of $V$
generated by $y^{4^2}$ as in the proof of Lemma \ref{ll4.4}. Let
${\mathcal P}^{4^2}$ and ${\mathcal Q}^{4^2}$ be the projections of
$V_{L}^{+}$  and $V$ to $M^{(4^2)}$ and $V^{(4^2)}$ respectively.
Let ${\mathcal I}^{4^2}(u^1,z)v^1={\mathcal P^{4^2}} Y(u^1,z)v^1$
for $u^1,v^1\in M^{(4)}$ be the intertwining operator of type
$$\left(\hspace{-3 pt}\begin{array}{c}
M^{(4^2)}\\
M^{(4)}\, \ \  M^{(4)}\end{array}\hspace{-3 pt}\right),$$ and
${\mathcal I}^{4^2}(\sigma_{1}(u^1),z)\sigma_{1}(v^1)={\mathcal
Q^{4^2}}Y(\sigma_{1}(u^1),z)\sigma_{1}(v^1)$ for $u^1,v^1\in
M^{(4)}$ be the intertwining operator of type
$$\left(\hspace{-3 pt}\begin{array}{c}
V^{(4^2)}\\
V^{(4)}\, \ \  V^{(4)}\end{array}\hspace{-3 pt}\right).$$ By the
proof of Lemma \ref{ll4.4}, the $L(1,0)$-module isomorphism
$\sigma_{1}$: $M^{(4)}\rightarrow V^{(4)}$ can be extended to the
$L(1,0)$-module isomorphism $\sigma_{1}$:
$$L(1,0)\oplus M^{(4)}\oplus M^{(4^2)} \rightarrow L(1,0)\oplus V^{(4)}\oplus
V^{(4^2)}$$ such that
\begin{equation}\label{eq4.11}\sigma_{1}(Y(u^1,z)v^1)=Y(\sigma_{1}(u^1),z)\sigma_{1}(v^1),
\end{equation}
\begin{equation}\label{eq4.12}\sigma_{2}(Y(u,z)w)=Y(\sigma_{1}(u),z)\sigma_{2}(w),
\end{equation}
for  $u^1,v^1\in M^{(4)}$, $u\in L(1,0)\oplus M^{(4)}\oplus
M^{(4^2)}$, $w\in M^{(k)}$.

It is well known that $M(1)^+=\oplus_{i\geq 0}L(1,(2i)^2)$ as
modules for the Virasoro algebra. We denote the $L(1,0)$-submodule
$L(1,r^2)$ of $M(1)^+$ by $M^{(r^2)}.$ Repeating the above procedure
we deduce that $U$ contains $L(1,0)$-submodules $V^{(r^2)}$
isomorphic to $L(1,r^2)$ for even $r.$  Then studying
$M^{(r^2)}\cdot M^{(s^2)}$ and $V^{(r^2)}\cdot V^{(s^2)}$, for
$(r,s)$ such that $r,s\in 2\Z_{+}$, $r\geq 2, s\geq 4$ and $r\leq
s,$ respectively (from small numbers to large ones), we conclude
that
$$
U=L(1,0)\bigoplus(\bigoplus_{r=1}^{\infty}L(1,(2r)^{2}),
$$
and there exists an $L(1,0)$-module isomorphism  $\sigma_{1}$ from
$M(1)^{+}$ to $U$ such that
\begin{equation}\label{eq4.13}\sigma_{1}(Y(u^1,z)v^1)=Y(\sigma_{1}(u^1),z)\sigma_{1}(v^1),
\end{equation}
\begin{equation}\label{eq4.14}\sigma_{2}(Y(u^1,z)w)=Y(\sigma_{1}(u^1),z)\sigma_{2}(w),
\end{equation}
for  $u^1,v^1\in M(1)^{+}$, $w\in M^{(k)}$. Then we prove that
\begin{theorem}\label{tt4.1}
$M(1)^{+}$ is isomorphic to the vertex operator subalgebra $U$
generated by $\omega$ and $J$.
\end{theorem}

\begin{lem}\label{ll4.5}
In the decomposition of $V$ into direct sum of irreducible
$L(1,0)$-modules, the multiplicity of $L(1,(2r)^{2})$ is 1, and the
multiplicity of $L(1,(2r+1)^{2})$ is zero,  for all $r\in\Z, r\geq
1$.
\end{lem}

\pf Suppose the lemma is false.  Let  $r$ be the smallest positive
integer such that there is a primary vector  with weight $r^{2}$
satisfying $u\notin U\cong M(1)^{+}$. Then $r\geq 2.$ Let $X$ be the
$U$-submodule of $V$ generated by $\sum\limits_{0\leq s<r^2}V_s$.
Then $X=U\oplus X^1$ where $X^1$ is the $U$-submodule generated by
the primary vectors whose weights are less than $r^2$ and are not
perfect square. Since $U$ is a simple vertex operator algebra, the
restriction of the bilinear form $(,)$ to $U$ is nondegenerate.
Clearly, the restriction of the bilinear form to $X^1$ is also
nondegenerate and $(U,X^1)=0.$ As a result, the restriction of the
bilinear form to $X$ is nondegenerate.

Let $X^{\perp}$ be the orthogonal complement of $X$ in $V.$ Using
the invariant property of the bilinear form we see that $X^{\perp}$
is also a $U$-module. Then $X^{\perp}=\sum_{s\geq r^2}X_s^{\perp}$
with $X_{r^2}^{\perp}\ne 0$ and each non-zero vector in
$X_{r^2}^{\perp}$ is a primary vector. Note that $J_tu=0$ for $t>3.$
Since $J_3$ preserves $X_{r^2}^{\perp},$ there exists $u\in
X_{r^2}^{\perp}$ satisfying
\begin{equation}\label{f1}
J_{3}u=c u \end{equation} for some $0\neq c\in\C$. Let $Z$ be the
$L(1,0)$-submodule generated by $u.$ Then $Z$ is isomorphic to
$L(1,r^2).$  By Theorem \ref{co2.2}, we see that $V^{(4)}\cdot
Z\subset Z\oplus_{t\geq (r+1)^2}X^{\perp}_t.$ In particular, $J_iu$
for $i=0,...,3$ is linear combination of $u, L(-1)u,$ $L_{-2}u,
L_{-1}^{2}u$ and $L_{-3}u,L_{-2}L_{-1}u,L_{-1}^{3}u.$

 Let
$$A_{k}=\{v\in V_{k}| L(n)v=0, n\geq 1\}.$$
By Theorem \ref{co2.2}, $A_{k}$ is invariant under the actions of
 $J_{3}$ and $u_{r^2-1},$ and $u_{i}|_{A_{k}}=0$, for $i\geq r^2.$
 Note that
 $$
 J_{3}u_{r^2-1}-u_{r^2-1}J_{3}=\sum\limits_{j=0}^{3}(J_{j}u)_{r^2+2-j}.$$
 It follows from the proof of Lemma \ref{ll2.4}  that for $j\geq 0$,  $(J_{j}u)_{r^2+2-j}|_{A_{k}}\in
 \C u_{r^2-1}|_{A_{k}}$. We deduce that on $A_{k}$,
 $$
 [J_{3},u_{r^2-1}]\in\C u_{r^2-1}.$$
So the  Lie subalgebra of $gl(A_{k},\C)$ generated by $J_{3}$
 and $u_{r^2-1}$ is solvable. Since $A_{k}$ is finite-dimensional, by the well-known Lie theorem, there
 exists $0\neq F\in A_{k}$ such that
 \begin{equation}\label{f2}J_{3}F, u_{r^2-1}F\in \C F.
 \end{equation}
Let $N$ be the $L(1,0)$-submodule of $V$ generated by $F.$ Then
$V^{(4)}\cdot N=N$ and $Z\cdot N=N.$ In particular, $J_{3}F\neq 0$,
$u_{r^2-1}F\neq 0$.

Let $W$ be the $M(1)^{+}$-submodule of $V$ generated by $u$. Then
$W$  has lowest weight $r^{2}.$ Although it is not clear that $W$ is
an irreducible $M(1)^+$-module, $W$ has a unique irreducible
quotient isomorphic to the irreducible module $M(1,\sqrt{2r})$ (see
\cite{DN1} for the notation). Since $W$ is a completely reducible
$L(1,0)$-module, it follows from \cite{DG} that $W$ contains an
$L(1,0)$-submodule $T$ isomorphic to $L(1,s^2)$ for some even $s\geq
r.$ Assume that $w$ is a nonzero primary vector of $T.$ Then $w$ is
not an element of $U=M(1)^+.$

By (\ref{f2}), we have
$$M(1)^{+}\cdot N=N, \ W\cdot N=N.$$
Let $v\in M(1)^+$ be a primary vector of weight $s^{2}$ and $T^1$
the irreducible $L(1,0)$-submodule of $V$ generated by $v.$ Then
both $T^1$ and $T$ are isomorphic to $L(1,s^{2}).$ By the above
discussion,
$$ T^1\cdot N=N, \ T\cdot N=N.$$
Similar to the proof of Lemma \ref{ll4.1}, we deduce that $v=c w$,
for some nonzero $c\in\C$, contradicting  the fact that $w\notin
M(1)^{+}$, $v\in M(1)^{+}$. So the lemma holds. \qed

\section{Characterization  of $V_{\Z\al}^{+}$}
\def\theequation{5.\arabic{equation}}
\setcounter{equation}{0}

Let $V$ be the vertex operator algebra as in Section 4. In this
section, we prove the main result of the paper. That is, $V$ is
isomorphic to $V_{\Z\alpha}^{+}$ with $(\al,\al)=2k$.

\vskip 0.3cm

Let $n$ be a positive integer. Set $X^1=e^{\beta}$, $X^2=e^{-\beta}$
where $(\beta,\beta)=2n.$ Let $W^{1}$ and $W^{2}$ be the irreducible
$M(1)^{+}$-submodules of $V_{\Z\beta}$ generated by $X^1$ and $X^2$
respectively. Then
\begin{equation}\label{fe1}
(X^{1},X^{1})=(X^2,X^2)=0, \ (X^1,X^2)=1
\end{equation}
 and
$$W^{2}\cdot W^{2}=W^{3}, \ W^{1}\cdot W^{3}=W^{2},$$
where $W^{3}$ is the irreducible $M(1)^{+}$-module generated by
$X^3=e^{-2\beta}$. Using the vertex operator algebra structure of
$V_{\Z\beta},$ we see that $Y(u^2,z)v^2$ for $u^2,v^2\in W^2$ and
$Y(u^1,z)v^3$ for $u^1\in W^1, v^3\in W^3$ are the intertwining
operators of type
$$\left(\hspace{-3 pt}\begin{array}{c}
W^3\\
W^2\, \ \  W^{2}\end{array}\hspace{-3 pt}\right),$$ and
$$\left(\hspace{-3 pt}\begin{array}{c}
W^2\\
W^1\, \ \  W^3\end{array}\hspace{-3 pt}\right)$$ for the vertex
operator algebra $M(1)^+.$ In particular, we have
$$X^2_{-2n-1}X^2=X^3, \ X^1_{4n-1}X^3=X^2,\ X^{2}_{i}X^{2}=0, \ \ i\geq -2n.$$
Using the explicit definition of $Y(X^1,z),$ Proposition 4.5.8 of
\cite{LL} and the fact that $X^2_mX^1=0$ for $m\geq 2n$, we have
\begin{align}
& (X^1_{2n-2}X^2)_{0}X^2\label{fl1}\\& =\sum\limits_{i=0}^{2n-2}(-1)^{i}\left(\begin{array}{c}2n-2\\i\end{array}\right)
(X^1_{2n-2-i}X^2_{i}-X^2_{2n-2-i}X^1_{i})X^2\nonumber\\
&
=\sum\limits_{i=0}^{2n-2}\sum\limits_{j=0}^{2n-i}\sum\limits_{s=0}^{2n}(-1)^{i+j+s+1}\left(\begin{array}{c}2n-2\\i\end{array}\right)
\left(\begin{array}{c}-2n\\j\end{array}\right)\left(\begin{array}{c}2n\\s\end{array}\right)X^1_{i+j+s}X^2_{2n-2-i-j-s}X^2\nonumber\\
& =2n X^2\neq 0.\nonumber
\end{align}

For $n\in\Z_{+}$, denote
$$A_{n}=\{v\in V_{n} | L(m)v=0, \ m\geq 1\}.$$ Then $\dim A_{n}<\infty$.
 Let $S$ be the set of  all the
positive integers $n$ such that $n\geq k$, $n$ is not a perfect
square and $\dim A_{n}\neq 0$. We prove in Section 4 that
\begin{equation}
V=M(1)^{+}\bigoplus(\bigoplus_{n\in S}(\dim A_{n})L(1,n)),
\end{equation}
where we identify the vertex subalgebra $U$ generated by $J$ and
$\omega$ with $M(1)^{+}$. Recall that the restriction of the
bilinear form to $M(1)^+$ is non-degenerate and the orthogonal
complement $M$ of $M(1)^{+}$ in $V$  is also an $M(1)^+$-module such
that $M\cap M(1)^+=0.$ Then clearly, as an $L(1,0)$-module,
$$M=\bigoplus_{n\in S}(\dim A_{n})L(1,n).$$

For $n\in S$,  by Theorem \ref{co2.2},
$$J_{3}A_{n}=A_{n}, \ J_{m}A_{n}=0, \ {\rm for} \ m\geq 4.$$
The following lemma is obvious.
\begin{lem}\label{ffl1}
Let $n\in S$ and $u\in A_{n}$ be such that $J_{3}u\in\C u$, then the
irreducible $L(1,0)$-module $L(1,n)$ generated by $u$ is an
irreducible $M(1)^{+}$-module.
\end{lem}

Note that the bilinear form on $V$ restricted to $A_{n}$ is still
non-degenerate.  We decompose $A_n$ into direct sum of
indecomposable $\C[J_3]$-modules: $A_n=\bigoplus_{i=1}^sX_i.$ Since
the form $(,)$ on $A_n$ is nondegenerate and $(J_3u,v)=(u,J_3v)$ for
$u,v\in A_n$ we may assume that
\begin{equation}\label{eq4.42}
(X_i,X_j)=0
\end{equation}
 if
$i\ne j.$

\begin{lem}\label{l5.2}
$J_{3}$ is semisimple on $A_{n}$.
\end{lem}

\pf  It is equivalent to show that each $X_i$ has dimension 1.
Suppose that $J_{3}$ is not semisimple on $A_{n}$. We may assume
that $\dim X_1>1$. Let
$$\{x^1,\cdots,x^r\}$$
be a basis of $X_1$ such that \begin{equation}\label{eq4.41}
(J_{3}-\la_{n}{\rm id})x^j=x^{j-1}
\end{equation}
 for
$j=1,2,\cdots,r$ where $x^0=0$ and $\la_n=4n^2-n$ (cf. \cite{DN1}).
Then
$$(J_{3}x^j,x^1)=(\la_{n} x^j+x^{j-1},x^1),$$
$$(x^{j},J_{3}x^1)=(x^j,\la_{n} x^1).$$
Since $(J_{3}x^j,x^1)=(x^j,J_{3}x^1)$, it follows that
\begin{equation}\label{eq5.2}
(x^{j-1},x^1)=0, \ j=2,\cdots,r.
\end{equation}
Using the non-degeneracy  of $(\cdot,\cdot)$,  we
 may choose $x^r$ such  that
\begin{equation}\label{eq5.2c}(x^r,x^1)=1,\
(x^r,x^r)=0.
\end{equation}
In fact, we may replace $x^r$ by $x^r-\frac{1}{2}(x^r,x^r)x^1$ if
$(x^r,x^r)\ne 0.$ Let $M^j$ be the irreducible $L(1,0)$-module
generated by $x^j$. Then  $M^1$ is an irreducible $M(1)^{+}$-module.
Using the invariant property of the bilinear form and (\ref{eq5.2})
we see that $(M^1,M^1)=0.$

{\bf Claim:} $P=M^1\cdot M^1=L(1,4n)$ is an irreducible
$M(1)^+$-module.

Note that $$ M^1\cdot M^1={\rm span}\{u_{l}v | u,v\in M^{1}, l\in{\mathbb
Z}\}$$ is an $M(1)^{+}$-module as $M^1$ is an $M(1)^+$-module. For
any $u, v\in M^1, w\in M(1)^{+}$ we have
$$(Y(u,z)v,w)=(v,Y(e^{zL(1)}(-z)^{(L(0)}v,z^{-1})w)=0$$
as the coefficients of $z^n$ in $Y(e^{zL(1)}(-z)^{(L(0)}v,z^{-1})w$
lie in $M^1.$ This implies that
$$(P,M(1)^+)=0$$
 and $P$ is an $M(1)^+$-submodule of $M.$ In particular,
$P$ is a direct sum of irreducible $L(1,0)$-modules whose lowest
weights are not perfect squares. Let $p\in\Z$ such that
$$x^1_px^1\ne 0, \ x^1_mx^1=0, \ m>p.$$
By Lemma \ref{ll2.4} we see that $x^1_px^1$ is the unique highest
weight vector up to a constant for the Virasoro algebra with the
highest weight $2n-p-1.$
 As a result, $x^1_px^1$ generates a highest weight module  for the Virasoro algebra and isomorphic to $L(1,2n-p-1)$ which is also an irreducible $M(1)^+$-module.
 Thus by Theorem \ref{co2.2} we have $M(1)^+$-module decomposition $P=L(1,2n-p-1)\oplus N$ where $N$ is the sum of highest weight
  modules for the Virasoro algebra whose highest weights
 are greater than $2n-p-1.$
Using the fact that  $(x^1,x^1)=0$ and Theorem \ref{tt33} we see
that $2n-p-1=4n.$

Let $N_m$ be the sum of irreducible highest weight modules for the Virasoro algebra in $N$ isomorphic to $L(1,m)$ for $m>4n.$
 Then $N_m$ is also an $M(1)^+$-module. If $N_m\ne 0$ for some $m>4n,$ then $N_m$ has a finite composition series as an $M(1)^+$-module. In particular,
  $N_m$ has a maximal $M(1)^+$-submodule such that the quotient is isomorphic to $M(1,\sqrt{2m}).$
   This  yields an intertwining operator of type $$\left(\hspace{-3 pt}\begin{array}{c} M(1,\sqrt{2m})\\
M(1,\sqrt{2n})\,M(1,\sqrt{2n})\end{array}\hspace{-3 pt}\right)$$ for
$m\ne 4n.$ This is a contradiction by Theorem \ref{tt33}. This
implies that $P=L(1,4n)$ is an irreducible $M(1)^+$-module.

Notice that $x^1_ix^1=0$ for $i\geq 0$ and $V_{1}=0$. We have
\begin{align*}
& (x^r_{2n-2}x^1)_{0}x^1\\
= &
\sum\limits_{i=0}^{2n-2}(-1)^{i+1}\left(\begin{array}{c}2n-2\\i\end{array}\right)
x^1_{2n-2-i}x^r_{i}x^1\\
=&
\sum\limits_{i=0}^{2n-2}\sum\limits_{j=0}^{2n-i}\sum\limits_{s=0}^{2n}(-1)^{i+j+s+1}\left(\begin{array}{c}2n-2\\i\end{array}\right)
\left(\begin{array}{c}-2n\\j\end{array}\right)\left(\begin{array}{c}2n\\s\end{array}\right)x^r_{i+j+s}x^1_{2n-2-i-j-s}x^1\\
= & 0.
\end{align*}
On the other hand, as $L(1,0)$-modules, $M^r\cong W^{1}$, $M^1\cong
W^{2}$, $P\cong W^{3}$ by module isomorphisms
\begin{align*} \tau_{1}: & M^r\rightarrow W^{1}, \
x^r\mapsto X^{1};\\
\tau_{2}: & {M}^1\rightarrow W^{2}, \ x^1\mapsto X^{2};\\
\tau_{3}: & P\rightarrow W^{3}, \ x^1_{-2n-1}x^1\mapsto X^{3}.
\end{align*}

Clearly, $M^r\cdot P$ is a direct sum of irreducible
$L(1,0)$-modules. It is easy to see from \cite{FZ} and \cite{L2}
that
$$
\dim_{ L(1,0)} \left(\hspace{-3 pt}\begin{array}{c} W^{2}\\
W^{3}\, \ W^{1}\end{array}\hspace{-3 pt}\right)=1.$$ It follows from
(\ref{fl1}) that
$$(x^r_{2n-2}x^1)_{0}x^1\neq 0$$
unless ${\cal P}(Y(u,z)w)=0$, for $u\in M^r$ and $w\in P$, where
${\cal P}(Y(u,z)w)$ is a projection of $Y(u,z)w$ to $M^1$. So
 $${\cal P}(Y(u,z)w)=0,$$ for $u\in M^r$ and
 $w\in P$. Then we have
\begin{align*}
& (x^r_{2n-1}x^1)_{-1}x^1\\
= &
\sum\limits_{i=0}^{2n-1}(-1)^{i}\left(\begin{array}{c}2n-1\\i\end{array}\right)
x^1_{2n-2-i}x^r_{i}x^{1}\\
= &
\sum\limits_{i=0}^{2n-1}\sum\limits_{j=0}^{2n-i}\sum\limits_{s=0}^{2n}(-1)^{i+j+s}
\left(\begin{array}{c}2n-1\\i\end{array}\right)
\left(\begin{array}{c}-2n\\j\end{array}\right)\left(\begin{array}{c}2n\\s\end{array}\right)x^r_{i+j+s}x^1_{2n-2-i-j-s}x^1\\
 = & 0.
\end{align*}
But by (\ref{eq5.2c}), we have
$$
(x^r_{2n-1}x^1)_{-1}x^1=x^1\neq 0,$$ a contradiction. The proof is
complete. \qed

By Lemmas \ref{ffl1} and \ref{l5.2}, we immediately have
\begin{lem}\label{fjl3}
As an $M(1)^{+}$-module, $V$ is completely reducible.
\end{lem}

It follows from Lemma \ref{l5.2} that there is a primary vector
$F\in V_{k}$ satisfying $J_{3}F=(4k^2-k)F$ and $(F,F)=2.$ Let
$V^{0}$ be the vertex operator subalgebra of $V$ generated by
$\omega$, $J$ and $F$.

\begin{lem}\label{ll4.9}
$V^{0}$ is linearly spanned by
$$L(-m_{1})\cdots L(-m_{s})J_{-n_{1}}\cdots
J_{-n_{t}}{\bf 1}, \ \ L(-m_{1})\cdots L(-m_{s})J_{-n_{1}}\cdots
J_{-n_{t}}F_{-p_{1}}\cdots F_{-p_{l}}F$$ where $m_{1}\geq m_{2}\geq
\cdots m_{s}\geq 1, n_{1}\geq n_{2}\geq \cdots n_{t}\geq 1,$ and
$p_{1}\geq p_{2}\geq \cdots p_{l}\geq 1$.
\end{lem}
\pf Clearly, $V^0$ is spanned by
$$L(m_{1})\cdots L(m_{s})J_{n_{1}}\cdots
J_{n_{t}}F_{p_{1}}\cdots F_{p_{l}}F$$ where
$m_{i},n_{j},p_{l}\in\Z$.

Recall that $M(1)^+=U$ is generated by $\omega$ and $J.$
Furthermore, $J_{i}J\in L(1,0)+V^{(4)}$ and $J_{i}F\in V^{(k)}$ for
$i\geq 0$ where $V^{(4)}, V^{(k)}$ are the irreducible
$L(1,0)$-modules generated by $J$ and $F,$ repectively. Using
Theorem \ref{tt33} and Lemma \ref{fjl3} we see that $F_{i}F\in U$
for $i\geq 0.$ The Lemma then follows from the commutator relations
$$
[J_{m},J_{n}]=\sum\limits_{i=0}^{\infty}\left(\begin{array}{c}m\\i\end{array}\right)(J_{i}J)_{m+n-i}
,$$
$$[F_{m},F_{n}]=\sum\limits_{i=0}^{\infty}\left(\begin{array}{c}m\\i\end{array}\right)(F_{i}F)_{m+n-i}$$
and
$$
[J_{m},F_{n}]=\sum\limits_{i=0}^{\infty}\left(\begin{array}{c}m\\i\end{array}\right)(J_{i}F)_{m+n-i}
$$ immediately. \qed

For $r,s\in\Z_{+}$ such that $r\neq s$, let $X,Y$ be irreducible
$M(1)^+$-submodules isomorphic to $L(1,r^{2}k)$ and $L(1,s^{2}k),$
respectively. It follows from Theorem \ref{tt33} and Lemma
\ref{fjl3} that $X\cdot Y$ is an $M(1)^+$-submodule of $V$
isomorphic to either $L(1,(r-s)^{2}k), L(1,(r+s)^{2}k)$ or
$L(1,(r-s)^{2}k)\bigoplus L(1,(r+s)^{2}k).$ Using Lemma \ref{ll4.9}
and Theorem \ref{tt33} gives
\begin{equation}\label{eq4.16}
V^{0}=(\bigoplus_{m=0}^{\infty}L(1,(2m)^{2})\bigoplus(\bigoplus_{r=1}^{\infty}c_{r}L(1,r^{2}k)),
\end{equation}
where $c_{r}\geq 0$.

\begin{lem}\label{ll4.15}
$c_{r}\leq 1.$
\end{lem}
\pf We prove the result by induction on $r$. By Lemma \ref{ll4.9},
$F$ is the only one linearly independent primary vector of weight
$k$ in $V^{0}$. That is, $c_1=1.$ Recall that $V^{(k)}$ is the
irreducible $M(1)^+$-module generated by $F$ and $V^{(k)}=L(1,k)$.
By Lemma \ref{ll2.4}, we have
$$
V^{(k)}\cdot V^{(k)}={\rm span}\{L(-m_{1})\cdots L(-m_{s})F_{n}F |
m_{1},\cdots, m_{s}\in {\mathbb Z}_{+}, n\in{\mathbb Z}\}.$$ Notice
that by Theorem \ref{tt33} and Lemma \ref{fjl3}, $F_{i}F\in U$ for
$i\geq -2k$. So $F_{m}F_{i}F\in V^{(k)}$ for $m\in\Z, i\geq -2k$.
From Lemma \ref{ll4.9} we see that any $u\in V^{0}_{4k}$ can be
written as
$$u=x+a F_{-2k-1}F$$
for some $x\in U+V^{(k)}, a\in \C$. Let
$$v^{1}=x^{1}+a_{1}F_{-2k-1}F, \ v^{2}=x^{2}+a_{2}F_{-2k-1}F$$ be two primary vectors of weight
$4k$ where $x^{1},x^{2}\in U\bigoplus V^{(k)}$, $a_{1},a_{2}\in \C$.
Since there is no primary vector of weight $4k$ in
$U+V^{(k)},$$a_{1},a_{2}$ are nonzero. Note that
$$a_{2}v^{1}-a_{1}v^{2}=a_{2}x^{1}-a_{1}x^{2}\in U\bigoplus V^{(k)}$$
is either zero or a primary element of weight $4k$. This implies
that $a^{2}v^{1}-a^{1}v^{2}=a^{2}x^{1}-a^{1}x^{2}=0$. That is, there
exists at most one linearly independent primary element of weight
$4k$.

By Theorem \ref{tt33} and Lemma \ref{fjl3}, for $p_{1}\in\Z_{+}$,
$$F_{-p_{1}}F=u^1+\sum\limits_{i=1}^{r^{(1)}_{1}}c^{(1)}_{i}L(-m^{(1)}_{i1})\cdots
 L(-m^{(1)}_{is_{i1}})F_{-2k-1}F,$$
for some $u^1\in U$, $c^{(1)}_{i}\in\C$, $m^{(1)}_{i1}\geq
m^{(1)}_{i2}\cdots \geq m^{(1)}_{is_{i1}}\geq 1.$ If there exists no
primary element of weight $4k$, then
$$V^{0}=U\bigoplus V^{(k)}.$$ The lemma holds.

Otherwise, let $v^{(4k)}$ be the primary element of weight $4k$
(unique up to a scalar). Then $v^{(4k)}$ generates an irreducible
$U$-module $V{(4k)}$ which is isomorphic to $L(1,4k)$ as
$L(1,0)$-modules.

Note that
$$F_{m}F_{-2k-1}F\in U+V^{(k)}+V^{(4k)}$$
 for $m\geq -4k.$
Then by Lemma \ref{ll4.9}, any element in $V^{0}$ of weight $9k$ is
a linear combination of an element in $U+V^{(k)}+V^{(4k)}$ and
$F_{-4k-1}F_{-2k-1}F$. Now let $v^{(9k)}$ be  a primary element of
weight $9k$, then we may assume that
$$v^{(9k)}=u+aF_{-4k-1}F_{-2k-1}F$$
for some $u\in U\bigoplus V^{(k)}\bigoplus V^{(4k)}$ and $0\neq
a\in\C$. As above we can prove that there exists at most one
linearly independent  primary element of weight $9k$ and for
$p_{1},p_{2}\in\Z_{+}$
$$F_{-p_{2}}F_{-p_{1}}F=\sum\limits_{i=1}^{t^{(2)}_{1}}d^{(2)}_{i}L(-n^{(2)}_{i1})\cdots
L(-n^{(2)}_{it_{i2}})F$$
$$\ \ \ \ \ \ \ \ \
+\sum\limits_{i=1}^{r^{(2)}_{1}}c^{(2)}_{i}L(-m^{(2)}_{i1})\cdots
L(-m^{(2)}_{is_{i2}})F_{-4k-1}F_{-2k-1}F$$ for some $c^{(2)}_{i},
d^{(2)}_{j}\in\C$, $m^{(2)}_{i1}\geq m^{(2)}_{i2}\cdots \geq
m^{(2)}_{is_{i2}}\geq 1, n^{(2)}_{i1}\geq n^{(2)}_{i2}\cdots \geq
n^{(2)}_{it_{i2}}\geq 1.$

Now assume that there exists only one linearly independent primary
vector $v^{(m^2k)}$ in $V^{0}$  of weight $m^{2}k$, for some $m\geq
3$, and for $p_{m-1},\cdots,p_{1}\in\Z_{+}$,
\begin{eqnarray*}
& F_{-p_{m-1}}\cdots F_{-p_{2}}F_{-p_{1}}F\\
& =
u+\sum\limits_{i=1}^{m-1}\sum\limits_{j=1}^{l_{ij}}c_{ij}L(-m^{(i)}_{j1})\cdots
L(-m^{(i)}_{jl_{ij}})F_{-2ik-1}\cdots F_{-4k-1}F_{-2k-1}F,
\end{eqnarray*}
where $u\in U\bigoplus V^{(k)}$, $c_{ij}\in \C$, $m^{(i)}_{j1}\geq
m^{(i)}_{j2}\geq \cdots\geq m^{(i)}_{jl_{ij}}\geq 1$.  Then by
Theorem \ref{co2.2}, $J_{3}v^{(m^2k)}\in\C v^{(m^2k)}$. So
$v^{(m^2k)}$ generates an irreducible $U$-module $V^{(m^2k)}$ which
is isomorphic to $L(1,m^2k)$ as $L(1,0)$-modules. By Theorem
\ref{tt33},
$$
F_{n}F_{-2(m-1)k-1}\cdots F_{-4k-1}F_{-2k-1}F\in
U\bigoplus(\bigoplus_{i=1}^{m}V^{(i^{2}k)}) \ {\rm for} \ n\geq
-2mk.$$ So any homogeneous element in $V^{0}$ of weight $(m+1)^{2}k$
is a linear combination of an element in
$U\bigoplus(\bigoplus_{i=1}^{m}V^{(i^{2}k)})$ and $F_{-2mk-1}\cdots
F_{-4k-1}F_{-2k-1}F$. \\ Let
$$v^{1}=x^{1}+a_{1}F_{-2mk-1}\cdots
F_{-4k-1}F_{-2k-1}F, \ v^{2}=x^{2}+a_{2}F_{-2mk-1}\cdots
F_{-4k-1}F_{-2k-1}F$$ be two primary vectors of weight $(m+1)^{2}k$,
 where $x^{1},x^{2}\in U\bigoplus(\bigoplus_{i=1}^{m}V^{(i^{2}k)})$, $a_{1},a_{2}\in
\C$. If $a_{1}=a_{2}=0$, then $x^{1}=x^{2}=0$. So we may assume that
$a_{1}\neq 0,a_{2}\neq 0$. Note that
$$a_{2}v^{1}-a_{1}v^{2}=a_{2}x^{1}-a_{1}x^{2}\in U\bigoplus(\bigoplus_{i=1}^{m}V^{(i^{2}k)})$$
is either zero or a primary vector of weight $(m+1)^{2}k$. This
implies that $a^{2}v^{1}-a^{1}v^{2}=a^{2}x^{1}-a^{1}x^{2}=0$. This
proves that there exists at most one linearly independent primary
vector of weight $(m+1)^{2}k$.

\qed

\begin{remark}\label{r5.6}
We actually also prove in Lemma \ref{ll4.15} that $V^{0}$ is
linearly spanned by
$$L(-m_{1})\cdots L(-m_{s})J_{-n_{1}}\cdots
J_{-n_{t}}{\bf 1}, \  \ L(-m_{1})\cdots L(-m_{s})F,$$ $$\
L(-m_{1})\cdots L(-m_{s})F_{-2mk-1}\cdots F_{-2k-1}F,$$
 where
$m_{1}\geq m_{2}\geq \cdots m_{s}\geq 1$, $m\geq 1$, $s,t\geq 0$.
\end{remark}

\begin{lem}\label{ll4.24} Let $V_{L}$ be the vertex operator algebra
 associated to the positive definite even lattice $L=\Z \al$ such that
 $(\a, \al)=2k$. Then
$V^0\cong V_{L}^{+}$.
\end{lem}
\pf For $m\in\Z_{+}$, set
$$E^{m}=e^{m\al}+e^{-m\al},\ \ E=E^1$$
and denote $N^{(m^2k)}\cong M(1,m\sqrt{2k})\cong L(1,m^{2}k)$ the
irreducible $M(1)^{+}$-module generated by $E^{m}$. Then
$$V_{L}^{+}=M(1)^{+}\bigoplus(\bigoplus_{m\in\Z_{+}}N^{(m^2k)}).$$
It is known that $(E^{m},E^{m})=2$ and
$J_{3}E^{m}=(4m^{4}k^2-m^{2}k)E^{m}$.

 By (\ref{eq4.16}) and Lemma \ref{ll4.15}, we may assume that
$$V^0=M(1)^{+}\bigoplus(\bigoplus_{m\in\Z_{+}}c_{m}V^{(m^{2}k)}),$$
where  $c_{m}\leq 1$, $V^{(m^{2}k)}\cong M(1,m\sqrt{2k})$. Since
$(E,E)=(F,F)=2$ we may assume that
$$
E_{-2k-1}E=u+aE^{2}, $$ $$ F_{-2k-1}F=u+v,$$ where $u\in M(1)^{+}$,
$0\ne a\in\C$ (see Lemma \ref{ll4.15}) and $v\in V$ is either zero
or a primary vector of weight $4k.$  Note that
$$J_{3}E=(4k^{2}-k)E, \ J_{3}F=(4k^2-k)F,$$
$$
(E_{-2k-1}E,E_{-2k-1}E)=(E,\sum\limits_{j=0}^{4k-1}\left(\begin{array}{c}4k-1\\j\end{array}\right)
(E_{j}E)_{2k-2-j}E),$$$$
(F_{-2k-1}F,F_{-2k-1}F)=(F,\sum\limits_{j=0}^{4k-1}\left(\begin{array}{c}4k-1\\j\end{array}\right)
(F_{j}F)_{2k-2-j}F).$$ Since $E_{j}E=F_{j}F\in M(1)^{+}$ for $j\geq
0$, it follows that
$$
(E_{-2k-1}E,E_{-2k-1}E)=(F_{-2k-1}F,F_{-2k-1}F).$$ So $$
(u+aE^{2},u+aE^{2})=(u+v,u+v)=(u,u)+2a^{2}=(u,u)+(v,v).$$ This
proves that $v\ne 0$. So we have $v=bF^2$ where $0\ne b\in \C$ and
$F^2\in V^0_{4k}$ is a primary vector such that $(F^2,F^2)=2.$ It
follows that $c_{2}=1$. If $a=-b$, we replace $F^{2}$ by $-F^{2}$.
So we may assume that $a=b$.

Let $\sigma$ be an $M(1)^{+}$-module isomorphism from
$$P^{(2)}=M(1)^{+}\bigoplus N^{(k)}\bigoplus N^{(4k)} \ \  to \ \
Q^{(2)}=M(1)^{+}\bigoplus V^{(k)}\bigoplus V^{(4k)}$$ defined by
$$
\sigma(E^{i})=F^{i}, \sigma(u)=u, \ i=1,2$$ for $u\in M(1)^+$ where
$F^1=F.$  For $u,v\in M(1)^{+}\bigoplus N^{(k)}$ and $n\in\Z$, we
have from Theorem \ref{tt33} that
$$
\sigma(u_{n}v)=(\sigma u)_{n}(\sigma v).$$

Now assume that $c_{i}=1$, for $i=1,2,\cdots,r$ and there exist
irreducible $M(1)^{+}$-modules $V^{(i^2k)}$ generated by non-zero
elements $F^{i}\in A_{i^2k}, i=1,2,\cdots, r$ such that $F^{1}=F$,
$(F^{i},F^{i})=2$ and the $\sigma$ can be extended to an
$M(1)^{+}$-module isomorphism from
$$P^{(r)}=M(1)^{+}\bigoplus (\bigoplus_{i=1}^{r}N^{(i^2k)}) \ \ to \
\ Q^{(r)}=M(1)^{+}\bigoplus (\bigoplus_{i=1}^{r}V^{(i^2k)})$$
satisfying
$$\sigma(E^{i})=F^{i}, \ i=1,2,\cdots,r$$
and for $u,v\in P^{(r)}$, $n\in\Z$ with $u_{n}v\in P^{(r)}$,
$$
\sigma(u_{n}v)=(\sigma u)_{n}(\sigma v).$$ Note that
$E^{1}_{-2rk-1}E^{r}\in P^{(r)}+N^{(r+1)}$, $F^{1}_{-2rk-1}F^{r}\in
Q^{(r)}$ or $ Q^{(r)}+V^{((r+1)^2k)}$ where $V^{((r+1)^2k)}$ is an
irreducible $M(1)^+$-module generated by a primary vector $F^{r+1}$
of weight $k(r+1)^2.$ By inductive assumption, we may assume that
$$
E^{1}_{-2rk-1}E^{r}=x+a_{1}E^{r+1}, \
F^{1}_{-2rk-1}F^{r}=\sigma(x)+b_{1}F^{r+1},$$ where $x\in P^{(r)}$,
$0\neq a_{1}\in\C$, $b_{1}\in\C$. We have
\begin{eqnarray*}
&(E^{1}_{-2rk-1}E^{r},E^{1}_{-2rk-1}E^{r})=(E^{r},E^{1}_{2rk+2k-1}E^{1}_{-2rk-1}E^{r})\\
&
=(E^{r},\sum\limits_{j=0}^{\infty}\left(\begin{array}{c}2rk+2k-1\\j\end{array}\right)(E^{1}_{j}E^{1})_{2k-2-j}E^{r})+(E^{r},E^{1}_{-2k-1}E^{1}_{2rk+2k-1}E^{r})
\end{eqnarray*}
\begin{eqnarray*}
&(F^{1}_{-2rk-1}F^{r},F^{1}_{-2rk-1}F^{r})=(F^{r},F^{1}_{2rk+2k-1}F^{1}_{-2rk-1}F^{r})\\
&
=(F^{r},\sum\limits_{j=0}^{\infty}\left(\begin{array}{c}2rk+2k-1\\j\end{array}\right)(F^{1}_{j}F^{1})_{2k-2-j}F^{r})+(F^{r},F^{1}_{-2k-1}F^{1}_{2rk+2k-1}F^{r})
\end{eqnarray*}
Notice that $E^{1}_{j}E^{1}=F^{1}_{j}F^{1}$, for $j\geq 0$ and
$$E^{1}_{2rk+2k-1}E^{r}\in N^{((r-1)^2k)}, \ F^{1}_{2rk+2k-1}F^{r}\in
V^{((r-1)^2k)}$$ Then by the inductive assumption, we have
$$(E^{1}_{-2rk-1}E^{r},E^{1}_{-2rk-1}E^{r})=(F^{1}_{-2rk-1}F^{r},F^{1}_{-2rk-1}F^{r}).$$
This implies that
$$
a_{1}^{2}=b_{1}^{2}\neq 0.$$ Thus  $c_{r+1}\neq 0.$ Replacing
$F^{r+1}$ by $-F^{r+1}$ if necessary we may assume that
$a_{1}=b_{1}$.

We extend the $M(1)^{+}$-module $\sigma$ to an $M(1)^{+}$-module
isomorphism from
$$
P^{(r+1)}=M(1)^{+}\bigoplus (\bigoplus_{i=1}^{r+1}N^{(i^2k)}) \ \ to
\ \ Q^{(r+1)}=M(1)^{+}\bigoplus (\bigoplus_{i=1}^{r+1}V^{(i^2k)})$$
satisfying
$$
\sigma E^{i}=F^{i}, \ i=1,2,\cdots, r+1$$ and for $u\in N^{(1)}$,
$v\in P^{(r)}$ and $n\in \Z$,
\begin{equation}\label{ee8}
\sigma(u_{n}v)=(\sigma u)_{n}(\sigma v).\end{equation}
 Now consider
$E^{2}_{-4rk+4k-1}E^{r-1}\in P^{(r)}+N^{((r+1)^{2}k)}$. From the
above discussion we know that
$$
E^{2}=u+aE^{1}_{-2k-1}E^{1}, \ F^{2}=u+aF^{1}_{-2k-1}F^{1}, $$where
$u\in M(1)^{+}$, $0\neq a\in\C$.  Note that for any $n\in \Z$
$$E^2_n=u_n+\sum_{i,j\in\Z}a_{i,j}E_iE_j,\ F^2_n=u_n+\sum_{i,j\in\Z}a_{i,j}F_iF_j$$
for some  $a_{i,j}\in \C.$ Combining this with (\ref{ee8}) yields
that for $u\in N^{(4k)}$, $v\in N^{((r-1)^{2}k)}$ and $n\in\Z$,
\begin{equation}\label{ee9}\sigma(u_{n}v)=(\sigma u)_{n}(\sigma
v).\end{equation} Continuing in this way, we show that for $u,v\in
P^{(r+1)}$, $n\in\Z$ with $u_{n}v\in P^{(r+1)},$
$$\sigma(u_{n}v)=(\sigma u)_{n}(\sigma v).$$
As a result,
$$
V^0=M(1)^{+}\bigoplus(\bigoplus_{m\in\Z_{+}}V^{(m^{2}k)}),$$ where
$V^{(m^2k)}$ is an irreducible $M(1)^{+}$-module generated by
non-zero element $F^{m}\in A_{m^2k}$ such that $(F^{m},F^{m})=2.$
Moreover, there exists a linear map $\sigma$ from $V_{L}^{+}$ to $V$
such that for $u,v\in V_{L}^{+}$, $n\in\Z$,
$$
\sigma(u_{n}v)=(\sigma u)_{n}(\sigma v).$$ Thus $V^{0}$ is
isomorphic to $V_{L}^{+}$. \qed

Here is the main theorem of this paper:
\begin{theorem} Let $V$ be as before. Then
$$V\cong V_{L}^{+}.$$
\end{theorem}
\pf By Lemma \ref{ll4.24}, $$V^{0}\cong V_{L}^{+}.$$ Then $V$ is an
extension of $V^0$ and $V$ is a direct sum of irreducible
$V^{0}$-modules. By the representation theory of $V_{L}^{+}$,
$V_{L}^{+}$ has only two irreducible modules with integral lowest
weights $0$ and $1$ \cite{DN2}. The assumption that $V_{1}=0$ and
$\dim V_0=1$ forces $V=V_{L}^{+}$. \qed

\end{document}